\documentclass[11pt, a4paper]{article}

\usepackage{amssymb,amsmath,amsthm,tocloft}
\usepackage{mathtools}
\usepackage{fancyhdr}
\usepackage{comment,enumerate}
\usepackage[american]{babel}
\usepackage{enumitem,verbatim}
\usepackage{stackengine}
\usepackage[usenames,dvipsnames,table]{xcolor}
\usepackage{graphicx,tikz,caption,subcaption}
\usepackage{bbm}

\usepackage{hyperref}
\usepackage[capitalise]{cleveref}

\usetikzlibrary{decorations.pathreplacing}
\usetikzlibrary{decorations.pathmorphing}
\usetikzlibrary{calc,arrows,decorations.markings}
\usetikzlibrary{shapes}
\tikzset{snake it/.style={decorate, decoration=snake}}
\usepackage{bbm,dsfont}

\usepackage[margin=2.75cm]{geometry}

\newtheorem{theorem}{Theorem}[section]
\newtheorem{prop}[theorem]{Proposition}
\newtheorem{lemma}[theorem]{Lemma}
\newtheorem{conjecture}[theorem]{Conjecture}
\newtheorem{cor}[theorem]{Corollary}

\newcommand{\repeatlabel}{}
\newtheorem*{repeatlemma}{Lemma \repeatlabel}

\theoremstyle{definition}

\newtheorem{obs}[theorem]{Observation}

\newtheorem{question}[theorem]{Question}

\newtheorem{claim}[theorem]{Claim}

\newcommand{\ex}{\mathrm{ex}}

\newcommand*{\abs}[1]{\lvert#1\rvert}

\newcommand{\e}{\epsilon}

\newcommand{\G}{\Gamma}

\newcommand{\R}{\mathcal{R}}

\DeclareMathOperator{\sat}{sat}
\DeclareMathOperator{\ssat}{ssat}
\DeclareMathOperator{\wsat}{wsat}

\title{Rainbow saturation for complete graphs}

\author{Debsoumya Chakraborti\thanks{
Discrete Mathematics Group, Institute for Basic Science (IBS), South Korea.
E-mail: {\tt \{debsoumya, kevinhendrey, benlund\}@ibs.re.kr}.
Debsoumya Chakraborti, Kevin Hendrey, and Ben Lund were supported by the Institute for Basic Science (IBS-R029-C1).}
\and 
Kevin Hendrey\footnotemark[1]\and
Ben Lund\footnotemark[1]\and
Casey Tompkins\thanks{
R\'enyi Institute, Hungary.
E-mail: {\tt casey@renyi.hu}.~Casey Tompkins was supported by the National Research, Development and Innovation Office, NKFIH (K135800).}}

\begin{document}

\maketitle

\begin{abstract}
We call an edge-colored graph \textit{rainbow} if all of its edges receive distinct colors. An edge-colored graph $\Gamma$ is called \textit{$H$-rainbow saturated} if $\G$ does not contain a rainbow copy of~$H$ and adding an edge of any color to $\Gamma$ creates a rainbow copy of $H$. 
The \textit{rainbow saturation number} $\sat(n,\R(H))$ is the minimum number of edges in an $n$-vertex $H$-rainbow saturated graph.
Gir\~{a}o, Lewis, and Popielarz conjectured that $\sat(n,\mathcal{R}(K_r))=2(r-2)n+O(1)$ for fixed $r\geq 3$.
Disproving this conjecture, we establish that for every $r\ge 3$, there exists a constant $\alpha_r$ such that
\[
r + \Omega\left(r^{1/3}\right) \le \alpha_r \le r + r^{1/2} 
\qquad \text{and} \qquad 
\sat(n,\R(K_r)) = \alpha_r n + O(1).
\]
Recently, Behague, Johnston, Letzter, Morrison, and Ogden independently gave a slightly weaker upper bound which was sufficient to disprove the conjecture.
They also introduced the \textit{weak rainbow saturation number}, and asked whether this is equal to the rainbow saturation number of $K_r$, since the standard weak saturation number of complete graphs equals the standard saturation number.
Surprisingly, our lower bound separates the rainbow saturation number from the weak rainbow saturation number, answering this question in the negative.
The existence of the constant $\alpha_r$ resolves another of their questions in the affirmative for complete graphs.
Furthermore, we show that the conjecture of Gir\~{a}o, Lewis, and Popielarz is true if we have an additional assumption that the edge-colored $K_r$-rainbow saturated graph must be rainbow.   
As an ingredient of the proof, we study graphs which are $K_r$-saturated with respect to the operation of deleting one edge and adding two edges. 
\end{abstract}

\section{Introduction}
Throughout this paper, given graphs $G$ and $H$ we say that $G$ is $H$-free if no subgraph of $G$ is isomorphic to $H$. The following two `dual' problems regarding $H$-free graphs are well-studied in extremal combinatorics. The classical extremal problem asks for the maximum number, denoted by $\ex(n,H)$, of edges in an $n$-vertex $H$-free graph.
Tur\'an~\cite{T41} solved this problem completely when $H$ is a complete graph~$K_r$. The graph saturation problem asks for the minimum number of edges in a maximal\footnote{maximal with respect to the subgraph relationship.} $H$-free graph with fixed number of vertices.
The following three notions of graph saturation have appeared in the literature numerous times throughout the last half a century.

\begin{itemize}[leftmargin=*]
    \item A graph $G$ is called \textit{$H$-saturated} if $G$ is a maximal $H$-free graph, i.e., $G$ is $H$-free and adding any edge to $G$ creates a copy of $H$. The \textit{saturation number} $\sat(n,H)$ denotes the minimum number of edges in an $n$-vertex $H$-saturated graph.
    
    \item A graph $G$ is called \textit{$H$-semisaturated} if adding any edge to $G$ creates a new copy of $H$. The \textit{semisaturation number} $\ssat(n,H)$ denotes the minimum number of edges in an $n$-vertex $H$-semisaturated graph.
    
    \item A graph $G$ is called \textit{weakly $H$-saturated} if there is an ordering $e_1,e_2,\ldots,e_t$ of the non-edges of $G$ such that for every $i\in \{1,2,\ldots,t\}$, the graph $(V(G),E(G)\cup \{e_1,e_2,\ldots,e_i\})$ contains a copy of $H$ using the edge $e_i$. The \textit{weak-saturation number} $\wsat(n,H)$ denotes the minimum number of edges in an $n$-vertex weakly $H$-saturated graph.
\end{itemize}
Observe that for every graph $H$, the following holds: 
\begin{equation}\label{sat_inequality}
    \sat(n,H) \ge \ssat(n,H) \ge \wsat(n,H).
\end{equation} 
The first result on graph saturation was obtained by Zykov~\cite{Z49} and independently by Erd\H{o}s, Hajnal, and Moon~\cite{EHM64}.
They proved that the unique $n$-vertex $K_r$-saturated graph with $\sat(n,K_r)$ edges is $K_{r-2}+\overline{K_{n-r+2}}$, the join of $K_{r-2}$ and an independent set of size $n-r+2$.
It is also well-known (see~\cite{A85,B65,B67,L77}) that for complete graphs $K_r$, the inequalities in \cref{sat_inequality} in fact hold with equality. 

\begin{theorem}[Bollobas~\cite{B65}, Erd\H{o}s--Hajnal--Moon~\cite{EHM64}, Zykov~\cite{Z49}]\label{ehm}
For every $n\ge r-2$,  
\[
\sat(n,K_r) = \ssat(n,K_r) = \wsat(n,K_r) = (r-2)(n-r+2)+\binom{r-2}{2}.
\]
Moreover, the unique graph witnessing the saturation number of $K_r$ is  $K_{r-2}+\overline{K_{n-r+2}}$.
\end{theorem}

For more results and problems related to graph saturation, we refer the readers to the survey by Faudree, Faudree, and Schmitt \cite{FFS11}.

In the setting of edge-colored graphs, a \emph{rainbow} copy of a graph $H$ (denoted $\R(H)$) is a graph isomorphic to $H$ whose edges are all assigned distinct colors.
There is a long history of studying rainbow analogs of classical results in combinatorics in the setting of edge-colored graphs. To mention a few of them in the context of Tur\'an-type extremal problems, see, e.g., \cite{ADGMS20,CKLLS22,GJ21,KMSV07,KSSV04} and the references therein.
Graph saturation was first studied in the setting of edge-colorings with a bounded number of colors by Hanson and Toft~\cite{HT87}, who focused on monochromatic cliques. In the same setting, Barrus, Ferrara, Vandenbussche, and Wenger~\cite{BFVW17} studied the saturation problem for rainbow copies of a given graph $H$.
Several further papers studied this problem for complete graphs and a few other simple graphs, such as paths, see, e.g.,~\cite{BFVW17,CMT20,FJL20,K18}.

In this paper, we focus on rainbow saturation problems in the setting of edge-colored graphs with infinitely many allowed colors (for concreteness, we can take the set of colors to be the set of natural numbers).
An edge-colored graph $\Gamma$ is $\mathcal{R}(H)$-saturated if it contains no $\mathcal{R}(H)$ subgraph and, for every non-edge $e$ and every color $c$, adding $e$ to $\Gamma$ with color $c$ creates a copy of $\mathcal{R}(H)$.
This version of rainbow saturation was first considered by Gir\~{a}o, Lewis, and Popielarz~\cite{GLP20}, who defined the \emph{rainbow saturation number} $\sat(n,\R(H))$ to be the minimum number of edges in an edge-colored $n$-vertex $\mathcal{R}(H)$-saturated graph. Note that if an edge-colored graph is $\R(K_r)$-saturated, then the underlying graph (ignoring the edge-coloring) has the property that adding any edge creates a new copy of $K_r$. 
Thus, Theorem~\ref{ehm} implies that whenever $n\ge r-2$, we have that $\sat(n,\R(K_r))\ge \ssat(n,K_r)= (r-2)n -\binom{r-1}{2}$. Gir\~{a}o, Lewis, and Popielarz~\cite{GLP20} conjectured the following.

\begin{conjecture}[\cite{GLP20}] \label{conjecture}
For every $r\ge 3$, there exists a constant $C_r$ depending only on $r$ such that, for any $n\ge 2(r-2)$, we have \[\sat(n,\R(K_r)) = 2(r-2)n + C_r.\]
\end{conjecture}

This conjecture was speculated based on the following construction. Consider the join graph between an independent set of size $n-2(r-2)$ and $K_{2(r-2)}$ minus a perfect matching. 
Then, color all the edges with distinct colors. Disproving Conjecture~\ref{conjecture} (it is also recently disproved in \cite{BJLMO22}, however, we obtain a slightly better upper bound) and also improving upon the trivial lower bound obtained by applying Theorem~\ref{ehm}, we obtain the following.

\begin{theorem} \label{thm1}
For each $r\geq 3$, there exist non-negative integers $\alpha_r,\beta_r$, and $N$ such that 
\begin{itemize}
\item $(r-2) + (1/2)(r-2)^{1/3} - 1/2\le \alpha_r\le (r-2) + \sqrt{r-1}$,
\item $\beta_r\le \binom{\alpha_r}{2}$, and 
\item for all $n\ge N$, we have $\sat(n,\R(K_r)) = \alpha_r(n-\alpha_r) + \beta_r.$
\end{itemize}
\end{theorem}

We believe that the value of $\alpha_r$ is closer to the upper bound; see the concluding remarks for a precise conjecture on this. For fixed $r$, calculating $\alpha_r$ amounts to solving a finite combinatorial problem, which we describe and analyze in Section~\ref{sec:ProofOfTheorem}.
In particular, we prove the following easy lower bound on~$\alpha_r$.
\begin{prop}\label{remark}
For each $r\geq 3$, the integer $\alpha_r$ defined in \cref{thm1} satisfies that
\begin{enumerate}
    \item $\alpha_r\ge r-1$ and
    \item $\alpha_r\ge r$ for every $r \ge 5$.
\end{enumerate}
\end{prop}

For sufficiently large $n$, this proposition already allows us to determine $\sat(n,\R(K_r))$ exactly for $r\in \{3,4,5\}$ and to determine $\alpha_r$ exactly for $r\in \{6,7,8,9\}$ (see \cref{thm:smallvalues}).

Behague, Johnston, Letzter, Morrison, and Ogden~\cite{BJLMO22} asked whether, for each graph $H$ there exists a constant $c(H)$ such that $\sat(n,H)=c(H)n + o(n)$.
\cref{thm1} answers this question for complete graphs.  
We remark that the upper bound in Theorem~\ref{thm1} only disproves Conjecture~\ref{conjecture} for $r\ge 4$. It turns out that Conjecture~\ref{conjecture} is true for $r=3$, as shown in Theorem~\ref{thm:semi-saturation}.

Similar to ordinary saturation, we next define semisaturation and weak-saturation in the edge-colored setting. 
Define the \emph{rainbow semisaturation number} $\ssat(n,\R(H))$ to be the minimum number of edges in an edge-colored $n$-vertex \textit{$\mathcal{R}(H)$-semisaturated graph}, where an edge-colored graph $\Gamma$ is called $\mathcal{R}(H)$-semisaturated if adding any edge with any color to $G$ creates a copy of $\mathcal{R}(H)$. 
Note that there is a $\mathcal{R}(H)$-semisaturated edge-colored graph with underlying graph $G$ if and only if $\mathcal{R}(G)$ is $\mathcal{R}(H)$-semisaturated.
On the other hand, Behague, Johnston, Letzter, Morrison, and Ogden~\cite{BJLMO22} defined an edge-colored graph $\Gamma$ to be \textit{weakly $\R(H)$-saturated} if there exists an ordering of the non-edges $e_1,e_2,\ldots,e_t$ of $\G$ such that, for any list of $c_1,c_2,\ldots,c_t$ of distinct colors, the non-edges $e_i$ in color $c_i$ can be added to $\G$, one at a time, so that every added edge creates a new copy of $\R(H)$. 
Then, define $\wsat(n,\R(H))$ to be the minimum number of edges in an edge-colored $n$-vertex weakly $\R(H)$-saturated graph.

As remarked by Behague, Johnston, Letzter, Morrison, and Ogden~\cite{BJLMO22}, the requirement to have distinct colors $c_1,c_2,\ldots,c_t$ excludes the possibility that all added edges have the same color, in which case the previously added edges do not contribute to making new copies of $\R(H)$, and the problem trivially reduces to rainbow semisaturation. 

Similar to the case of ordinary saturation, we observe the following.

\begin{obs}\label{obs:comparison}
For every graph $H$, we have
\[\sat(n,\R(H)) \ge \ssat(n,\R(H)) \ge \wsat(n,\R(H)).\]
\end{obs}

By a similar argument as before, it is easy to see that $\ssat(n,\R(K_r)) \ge \ssat(n,K_r)$. We precisely determine $\ssat(n,\R(K_r))$ for all $r$ and sufficiently large $n$ in the following theorem.

\begin{theorem}\label{thm:semi-saturation}
For every $r\ge 3$, there exists $N$ such that for all $n\ge N$, we have 
\begin{equation*}
    \ssat(n,\R(K_r))=\begin{cases}
    \sat(n,\R(K_3)) = 2(n-2), & \text{if $r=3$}.\\
    (r-1)(n-r+1)+\binom{r-1}{2}, & \text{otherwise}.
    \end{cases}
\end{equation*}
Moreover, the equality is uniquely achieved by $\mathcal{R}(\overline{K_2}+\overline{K_{n-2}})$ if $r=3$ and by ${\mathcal{R}(K_{r-1}+\overline{K_{n+1-r}})}$ otherwise.
\end{theorem}

Next, answering Question~6.2 of \cite{BJLMO22} in the negative, we show that rainbow versions of the weak-saturation and semisaturation numbers of complete graphs of order at least $5$ are smaller than the rainbow saturation number by an additive linear factor. This behavior is different from the ordinary saturation numbers of complete graphs as mentioned in \cref{ehm}. As a simple corollary of \cref{thm1,remark,thm:semi-saturation,obs:comparison}, we obtain the following.
\begin{cor}\label{cor:comparison}
For every $r\ge 5$, there exists $N$ such that for all $n\ge N$, we have $$\sat(n,\R(K_r)) - n/2 \ge \ssat(n,\R(K_r)) \ge \wsat(n,\R(K_r)).$$
\end{cor}

Now, we turn our attention to proving a variant of Conjecture~\ref{conjecture} with an extra condition. 
Notably, our constructions for the upper bound in Theorem~\ref{thm1} are not rainbow.
In fact, we show that a modified version of \cref{conjecture} is correct, with the added hypothesis that no color appears more than once in the $\R(K_r)$-saturated graph.
For any graph~$H$, let $\sat(\R(K_n),\R(H))$ denote the minimum number of edges in an $n$-vertex graph $G$ such that $\R(G)$ is $\R(H)$-saturated, if such a graph exists; otherwise, define $\sat(\R(K_n),\R(H))$ to be infinity\footnote{As an example, we show that $\sat(\R(K_n), \R(K_{1,4}))=\infty$ for odd $n\geq 5$. Given an $n$-vertex graph $G$, observe that if $\R(G)$ is $\R(K_{1,4})$-free, then $G$ has no vertex of degree more than $3$, and that some vertex $v$ has degree less than $3$ since $n$ is odd.
When adding an edge to $\R(G)$ incident to $v$, at most one copy of $K_{1,4}$ is created, and the color of the new edge may be chosen so that this copy is not rainbow.}.
Note that $\sat(n,\R(H))\le \sat(\R(K_n),\R(H))$ for every graph $H$. We have the following result.

\begin{theorem} \label{thm:rainbow_construction}
For every $r\ge 3$, there exists $N$ such that for all $n\ge N$, we have \[\sat(\R(K_n),\R(K_r)) = 2(r-2)(n-r+1).\] 
Moreover, the equality is uniquely achieved by
$\R(\overline{(r-2)K_2}+\overline{K_{n+4-2r}})$, the rainbow the complete multipartite graph with $r-2$ parts of size $2$ and one part of size $n+4-2r$. 
\end{theorem}

To prove \cref{thm:semi-saturation,thm:rainbow_construction}, we consider a variant of graph saturation that does not involve colors. A graph $G$ is \textit{$(H,1)$-saturated} if $G$ is $H$-saturated, and it is not possible to remove an edge and add two new edges without creating a copy of $H$ (see the concluding remarks for a discussion on a natural generalization of this notion).
Note that if $G$ is an $n$-vertex $H$-saturated graph with the maximum possible number of edges, then it is $(H,1)$-saturated.
Let $\sat_1(n,H)$ denote the minimum number of edges in an $n$-vertex $(H,1)$-saturated graph.

We use the following variant of semisaturation to prove Theorem~\ref{thm:semi-saturation}. 
A graph $G$ is \textit{$(H,1)$-semisaturated} if $G$ is $H$-semisaturated, and it is not possible to remove an edge and add two new edges without creating a copy of $H$. 
Let $\ssat_1(n,H)$ denote the minimum number of edges in an $n$-vertex $(H,1)$-semisaturated graph. 
The next statement highlights the relationship of $(H,1)$-saturation with rainbow saturation.

\begin{prop}\label{prop:comparison_semi}
For all graphs $H$ and $G$, if $\R(G)$ is $\R(H)$-semisaturated, then $G$ is $(H,1)$-semisaturated.
In particular, if $\R(G)$ is $\R(H)$-saturated, then $G$ is $(H,1)$-saturated.
\end{prop}

In light of \cref{prop:comparison_semi}, our strategy to prove \cref{thm:semi-saturation}
is to simultaneously determine $\ssat_1(n,K_r)$ and $\ssat(n,\R(K_r))$, and our strategy to prove \cref{thm:rainbow_construction}
is to simultaneously determine $\sat_1(n,K_r)$ and $\sat(\R(K_n),\R(K_r))$.
Thus we also obtain the following results.

\begin{theorem} \label{thm2'}
For every $r\ge 3$, there exists $N$ such that for all $n\ge N$, we have
\begin{equation*}
    \ssat_1(n,\R(K_r))=\begin{cases}
    2(n-2), & \text{if $r=3$}.\\
    (r-1)(n-r+1)+\binom{r-1}{2}, & \text{otherwise}.
    \end{cases}
\end{equation*}
Moreover, the equality is uniquely achieved by $\overline{K_2}+\overline{K_{n-2}}$ if $r=3$, and by $K_{r-1}+ \overline{K_{n-r+1}}$ otherwise.
\end{theorem}

\begin{theorem} \label{thm2}
For every $r\ge 3$, there exists $N$ such that for all $n\ge N$, we have \[\sat_1(n,K_r) = 2(r-2)(n-r+1).\] 
Moreover, the equality is uniquely achieved by $\overline{(r-2)K_2}+\overline{K_{n+4-2r}}$. 
\end{theorem}

A major ingredient of all of our results is a structural lemma which we prove in the next section.
Aside from unifying our proofs, the advantage of this lemma is that it gives us a form of structural stability for constructions that are within a small number of edges of the optimal bounds we obtain.
In particular, for $\sat_1(n,K_r)$ (and also for $\sat(\R(K_n),\R(K_r))$), there is a unique construction within $n-o(n)$ edges of the $1$-saturation number. 
Conversely, for $\sat(n,\R(K_r))$, we show in Section~\ref{sec:ProofOfTheorem} that for $n$ sufficiently large in terms of $r$, there are $n$-vertex $\R(K_r)$-saturated edge-colored graphs with any desired number of edges between $\sat(n,\R(K_r))$ and $\binom{n}{2}$.

\medskip
\noindent{\bf Organization.}
In the next section, we prove the structural lemma that will be used to prove the main results of this paper. In \cref{sec:ProofOfTheorem}, we prove \cref{thm1}. We also determine $\sat(n,K_r)$ for small values of $r$, provide some explicit constructions, and prove some related results about $\R(K_r)$-saturated graphs.
In \cref{sec:MoreProofs}, we prove \cref{prop:comparison_semi} and Theorems~\ref{thm:semi-saturation}, \ref{thm:rainbow_construction}, \ref{thm2'}, and~\ref{thm2}, obtaining stability results for edge-minimal constructions.
In \cref{sec:Conclusion}, we finish with several concluding remarks, which include (i) a few open questions and conjectures, (ii) a natural extension of the function $\sat_1(n,H)$ along with some preliminary observations.

\medskip
\noindent{\bf Notation.}
For a positive integer $s$, we use the standard notation $[s]$ to denote the set $\{1,2,\ldots,s\}$.
For a graph $G$, we denote its vertex set by $V(G)$, its edge set by $E(G)$, and the maximum degree by $\Delta(G)$.
We denote by $\omega(G)$ the \emph{clique number} of $G$, which is the maximum order of a clique in $G$.
For a graph $G$ and a set $S\subseteq V(G)$, let $G[S]$ denote the subgraph induced by $S$ and $G-S$ denote the subgraph induced by $V(G)\setminus S$. For $s\in V(G)$, we will write $G-s$ instead of $G-\{s\}$. 
For a vertex $v\in V(G)$ we denote by $N_G(v)$ the set of neighbors of $v$ and by $N_G[v]$ the set $\{v\}\cup N_G(v)$, and for a set $S\subseteq V(G)$ we denote by $N_G[S]$ the set $\bigcup\{N_G[v]:v\in S\}$. We may omit the subscripts if the graph is clear from the context.
We use the same notation for edge-colored graphs, where a subgraph $\Lambda$ of an edge-colored graph $\Gamma$ inherits the restricted edge-coloring.

\section{Central lemma}

We start with a key lemma which we use to prove our main results (i.e., Theorems~\ref{thm1}, \ref{thm:semi-saturation}, \ref{thm:rainbow_construction}, \ref{thm2'}, and~\ref{thm2}) in this paper. 
Expanding on ideas used in \cite{CCH21,CL20}, this lemma determines the global structure of a graph with a linear number of edges based on some local forbidden substructures, which reduces all of the problems we consider to simpler problems.
Not only does this unify the proofs of all of our main results, but it also gives us a form of structural stability each time we apply it.
We expect that it may also be useful for providing lower bounds in other graph saturation problems.

\begin{lemma}\label{lem:general}
 Given positive integers $k$ and $\ell$ and a positive real number $\epsilon$, there is some positive integer $N$ such that given $n\geq N$, an $n$-vertex graph $G$ with at most $(k+1-\epsilon)n$ edges and family $\mathcal{F}$ of $k$-vertex graphs, at least one of the following holds for some subset $S\subseteq V(G)$ of size $k$:
    \begin{enumerate}
        \item there is some independent set $X$ in $G-S$ of size at least $\ell$ and some $v\in X$ such that for every $w\in X\setminus \{v\}$, we have that $N(v)\cap N(w)\subseteq S$ and $G[N(v)\cap N(w)]$ is not isomorphic to a graph in $\mathcal{F}$, or
        \item $S$ is complete to $V(G)\setminus S$ and $G[S]$ is isomorphic to a graph in $\mathcal{F}$.
    \end{enumerate}
\end{lemma}

When applying the above lemma, we will always show that the first condition is impossible 
and thus obtain the second condition, which gives a lower bound on the number of edges.
We remark that the above lemma with $\e=1$ is sufficient for our main results. However, allowing $\e<1$ gives us information about constructions that are close to optimal, and in some cases, gives a stability result.
We can interpret this lemma as giving a $(k+1)n-o(n)$ lower bound on the number of edges in an $n$-vertex graph with no such set $S$ of size $k$.
For integers $k\geq 1$, $\ell\geq 2$, and $n\geq 2(k+1)$
the graph $K_{k+1,n-(k+1)}$ has $(k+1)n-(k+1)^2$ edges and does not satisfy either condition (regardless of the choice of $\mathcal{F}$ or $S$), since every pair of non-adjacent vertices has more than $k$ common neighbors and no set of $k$ vertices is complete to its complement. 
Thus our $(k+1)n-o(n)$ lower bound is tight up to the $o(n)$ term.
It would be nice to replace $o(n)$ with a constant depending only on $k$ and $\ell$, and we suspect that this is possible.

\begin{proof}[Proof of \cref{lem:general}]
Given $N$ sufficiently large, we will assume that the first condition fails and prove that the second condition holds.

We begin by defining $V_{\textrm{big}}(G)$ to be the set of vertices $v$ with $d(v)\geq n^{1/4}$, and $V_{\textrm{small}}(G)=V(G)\setminus V_{\textrm{big}}(G)$. 
Let $G_k$ be the graph with $V(G_k):=\{v\in V_{\textrm{small}}(G):|N(v)\cap V_{\textrm{big}}(G)|\leq k\}$, such that vertices $v$ and $w$ are adjacent if the distance from $v$ to $w$ in $G[V_{\textrm{small}}(G)]$ is at most~$2$.

Since $G$ has at most $(k+1)n$ edges, we have $|V_{\textrm{big}}(G)|\leq 2(k+1)n^{3/4}$.
Note that by the definition of $V(G_k)$, we have $|E(G)|\geq (k+1)|V_{\textrm{small}}(G)\setminus V(G_k)|$, and so 
\[
|V(G_k)|\geq |V_{\textrm{small}}(G)| - (k+1-\epsilon)n/(k+1) = \epsilon n/(k+1)-|V_{\textrm{big}}(G)|.
\]
Note that $\Delta(G_k)< (n^{1/4})+n^{1/4}(n^{1/4}-1)= \sqrt{n}$. 
For each $v\in V(G_k)$, define $B_v:=V(G_k)$ if $G[N(v)\cap V_{\textrm{big}}(G)]$ is not isomorphic to a graph in $\mathcal{F}$, and otherwise define 
\[B_v:=\{v\}\cup \{w\in V(G_k):|N(w)\cap N(v)\cap V_{\textrm{big}}(G)|<k\}.\]
If $|B_v|> (\ell-1)(\sqrt{n}+1)$, then we can find an independent set $X=\{v_0,\dots , v_{\ell-1}\}$ in $G_k[B_v]$ by setting $v_0:=v$ and iteratively selecting $v_i$ for $i\in [\ell-1]$ to be a vertex in $B_v\setminus N_{G_k}[\{v_0,\dots,v_{i-1}\}]$. 
This is possible since $|N_{G_k}[\{v_0,\dots,v_{i-1}\}]|\leq i(\Delta(G_k)+1)<|B_v|$.
Since $X$ is independent in $G_k$, we have $N_G(v_0)\cap N_G(v_i)\subseteq N_G(v_0)\cap V_{\textrm{big}}(G)$ for all $i\in [\ell-1]$, where the subset relation is strict if $G[N_G(v_0)\cap N_G(v_i)]$ is isomorphic to a graph in $\mathcal{F}$. Since all graphs in $\mathcal{F}$ have $k$ vertices, the first condition of the lemma is satisfied (with $S=N(v)\cap V_{\textrm{big}}(G)$).

We may therefore assume that $|B_v|\leq (\ell-1)(\sqrt{n}+1)$ for each $v\in V(G_k)$, and hence (given that $n$ is sufficiently large) for any pair of vertices $v,v'\in V(G_k)$, there is some $w\in V(G_k)\setminus (B_{v}\cup B_{v'})$. 
Thus, 
\[
k=|N(v)\cap V_{\textrm{big}}(G)|= |N(v)\cap N(w)\cap V_{\textrm{big}}(G)|=|N(v')\cap N(w)\cap V_{\textrm{big}}(G)|=|N(v')\cap V_{\textrm{big}}(G)|,
\]
and $G[N(v)\cap V_{\textrm{big}}(G)]$ is isomorphic to a graph in $\mathcal{F}$.
It follows that there is a set $S\subseteq V_{\textrm{big}}(G)$ of size $k$ such that $V(G_k)$ is complete to $S$, and $G[S]$ is isomorphic to a graph in $\mathcal{F}$.
Note that by the definition of $G_k$, each vertex $w\in V(G_k)$ satisfies $N(w)\cap V_{\textrm{big}}(G)=S$.

Now suppose there is a vertex $v\in V_{\textrm{small}}(G)\setminus V(G_k)$ which is not adjacent to every vertex of $S$.
We construct an independent set $X=\{v_0,\dots, v_{\ell-1}\}$ in $G[V(G_k)\cup \{v\}]$ by setting $v_0:=v$ and for each $i\in [\ell-1]$, choosing an arbitrary vertex $v_i\in V(G_k)$ 
such that $v_i$ is at distance more than 2 from $v$ in $G[V_{\textrm{small}}]$. This works since the number of vertices at distance at most 2 from $v$ in $G[V_{\textrm{small}}]$ is less than $\sqrt{n} \le |V(G_k)| - l$. 
Thus, the first condition is satisfied with the constructed $X$.

We may therefore assume that $V_{\textrm{small}}$ is complete to $S$. In particular, we have
\[|E(G-S)|\leq (k+1-\epsilon)n-|V_{\textrm{small}}(G)|k\leq (1-\epsilon)n+o(n).\]
Since each component $C$ of $G-S$ satisfies $|E(G[C])|\geq |C|-1$, it follows that $G-S$ has at least $\ell$ components provided $n$ is sufficiently large. 
Now, if some vertex $v\in V(G)\setminus S$ is not complete to $S$, we can find the independent set $X$ satisfying the first condition by taking each vertex from a different component of $G-S$ (with $v$ as the specially selected vertex).
Hence, if the first condition does not hold, then $S$ is complete to $V(G)\setminus S$, and the second condition is satisfied.
\end{proof}

\section{Rainbow Complete Graphs --- Proof of Theorem~\ref{thm1}}\label{sec:ProofOfTheorem}
In this section, we aim to prove \cref{thm1}. In the next subsection, we apply \cref{lem:general} to find an alternative description of $\alpha_r$ in \cref{thm1} as promised in the introduction.
In the subsequent two subsections, using this alternative formulation, we provide upper and lower bounds on $\alpha_r$. 
We note that while we cannot precisely determine $\alpha_r$ for $r\geq 10$, our machinery still gives us a surprising amount of structural information about $\R(K_r)$-saturated graphs which have close to the minimum possible number of edges.
In particular, these graphs contain $\R(K_{\alpha_r,n-\alpha_r})$ as a subgraph (see, \cref{cor:structuralstability}).
On the other hand, we conclude this section by proving a non-stability result (see, \cref{thm:non_stability}), namely that for sufficiently large $n$ and any $m$ between $\sat(n,\R(K_r))$ and $\binom{n}{2}$, there is a $\R(K_r)$-saturated edge-colored graph with exactly $n$ vertices and $m$ edges.

\subsection{Reducing {\boldmath $\mathcal{R}(K_r)$}-saturation to an equivalent problem}

Given a positive integer $k$, let $\hat{\mathcal{F}}_k$ be the class of all edge-colored graphs $\G$ satisfying the following properties.
\begin{enumerate}
    \item \label{item:prop1}$\Gamma$ does not contain $\mathcal{R}(K_{k+1})$,
    \item \label{item:prop2} for each $v \in V(\Gamma)$, the edge-colored graph induced by $V(\Gamma) \setminus \{v\}$ contains $\mathcal{R}(K_k)$, and
    \item \label{item:prop3} for every color $c$, $\Gamma$ contains a $\mathcal{R}(K_k)$ that does not use $c$.
\end{enumerate}
Define $f(k)$ to be the smallest integer such that there exists an edge-colored $f(k)$-vertex graph in $\hat{\mathcal{F}}_k$.
Let $\Lambda_k$ be an $f(k)$-vertex $\mathcal{R}(K_{k+1})$-saturated edge-colored graph in $\hat{\mathcal{F}}_k$ with as few edges as possible (such a graph exists since the $f(k)$-vertex graph in $\hat{\mathcal{F}}_k$ with the most edges is $\mathcal{R}(K_{k+1})$-saturated). 
Let $g(k) = |E(\Lambda_k)|$, and let $g'(k)$ be the minimum number of edges of an $f(k)$-vertex graph in $\hat{\mathcal{F}}_k$.

The motivation for these definitions is the following construction.
For positive integers $r\geq 3$ and $n\geq f(r-2)+2$, we define $\Gamma_{r,n}$ to be the edge-colored graph obtained by taking the complete join of $\Lambda_{r-2}$ with an independent set $I$ of size $n-f(r-2)$, and assigning all edges incident with $I$ unique colors which do not appear anywhere else in the graph.
Now Property~\ref{item:prop1} guarantees that $\Gamma_{r,n}$ is $\R(K_r)$-free, Properties~\ref{item:prop2} and~$\ref{item:prop3}$ ensure that adding any edge to $I$ in any color will create a copy of $\R(K_r)$, and the fact that $\Lambda_{r-2}$ is $\R(K_{r-1})$ saturated and $|I|\geq 2$ ensures that adding any edge not incident to $I$ in any color will create a copy of $\R(K_r)$.
Thus $\Gamma_{r,n}$ is $\R(K_r)$-saturated.
The main result of this subsection is that this construction is either optimal or at most $g(r-2)-g'(r-2)-1$ edges away from being optimal.

Given $r\geq 3$, let $\mathcal{F}_{r-2}$ be the class of all $f(r-2)$-vertex edge-colored graphs $\Gamma$ such that some subgraph of $\G$ is in $\hat{\mathcal{F}}_{r-2}$.

\begin{lemma}\label{lem:rainbowlower}
    Given positive integers $n\geq r\geq 3$ and an $n$-vertex $\R(K_r)$-saturated edge-colored graph $\G$, there is no set $S\subseteq V(\G)$ of size $f(r-2)$ for which there is an independent set $X$ in $\G-S$ of size $2(f(r-2)\binom{f(r-2)}{2})^{f(r-2)}$ such that for some $v\in X$ and every $w\in X\setminus \{v\}$, we have that $N(v)\cap N(w)\subseteq S$ and $\G[N(v)\cap N(w)]\notin \mathcal{F}_{r-2}$.
\end{lemma}
\begin{proof}  
Suppose for contradiction that such an $S$, $X$, and $v$ exist.
Let $S'$ be a maximal subset of $S$ such that there is a set $B_{S'}\subseteq X$ of size at least $2(f(r-2)\binom{f(r-2)}{2})^{f(r-2)-|S'|}$ with the property that for all $u\in S'$ and $w,w'\in B_{S'}$, we have that $uw$ and $uw'$ are edges of $\G$ and have the same color.
If $S'=\emptyset$, we take $B_{S'}=X$.
Note that $|S'|\leq |S|= f(r-2)$ and so $|B_{S'}|\geq 2$.

Fix a vertex $v'\in B_{S'}$, with $v'=v$ if $S'=\emptyset$, and consider an arbitrary vertex $w\in B_{S'}\setminus \{v'\}$.
Consider the edge-colored graph $\Lambda_{v',w}$ obtained from $\G[(N(v')\cap N(w))\setminus S']$ by deleting every edge $uu'$ for which $uv'$ has the same color as $u'v'$ or $uw$ has the same color as $u'w$.
Now, if there is some vertex $u\in V(\Lambda_{v',w})$ such that $\Lambda_{v',w}-u$ is $\mathcal{R}(K_{r-2})$-free, then we can add $v'w$ to $\G$ and assign it the same color as $v'u$ without creating a copy of $\mathcal{R}(K_r)$, a contradiction. 
Similarly, if there is some color $c$ which is used by every copy of $\mathcal{R}(K_{r-2})$ in $\Lambda_{v',w}$, then we can add $v'w$ to $\G$ and assign it color $c$ without creating $\mathcal{R}(K_r)$, a contradiction.
Note that by our assumption and choice of $v'$, $\Lambda_{v',w}\notin \hat{\mathcal{F}}_{r}$, since if $v'\neq v$, then $|V(\Lambda_{v',w})|<f(r-2)$.
Hence, by the definition of $\hat{\mathcal{F}}_{r}$, we conclude that $\Lambda_{v',w}$ contains $\mathcal{R}(K_{r-1})$. 
Now given that $\G$ does not contain $\mathcal{R}(K_r)$, for each $w'\in \{v',w\}$ there is some edge $e\in \Lambda_{v',w}$ and some vertex $u\in V(\Lambda_{v',w})$ such that $w'u$ has the same color as $e$.
Hence, for each $w\in B_{S'}$, there is some $e\in \G[S]$ and some vertex $u\in S\setminus S'$ such that $wu$ has the same color as $e$.
By the pigeonhole principle, there is a vertex $u\in S\setminus S'$ and a subset $B_{S''}$ of $B_{S'}$ of size at least $|B_{S'}|/(f(r-2)\binom{f(r-2)}{2})$ such that for all $w,w'\in B_{S''}$, we have that $uw$ and $uw'$ are edges of $\G$ and have the same color. Now $S'':=S'\cup \{u\}$ contradicts the maximality of $S'$, which completes the proof.
\end{proof}

Now the following is an immediate corollary of \cref{lem:general,lem:rainbowlower}, where we take $\mathcal{F}$ to be the class of all underlying graphs of edge-colored graphs in $\mathcal{F}_{r-2}$. 
\begin{cor}\label{cor:structuralstability}
    For every integer $r\geq 3$ and real number $\e >0$, there is some positive integer $N$ such that for every $n\ge N$, every $n$-vertex edge-colored graph $\G$ with at most $(f(r-2)+1-\e)n$ edges which is  $\R(K_r)$-saturated has a set $S$ of $f(r-2)$ vertices such that $S$ is complete to $V(\G)\setminus S$ and $|E(\G[S])|\geq g'(r-2)$.
\end{cor}

We now prove the main result of this subsection, which shows that the constant $\alpha_r$ in \cref{thm1} is always equal to $f(r-2)$.
\begin{lemma}\label{lem:sat_reduction}
    For every integer $r\geq 3$, there exist integers $\beta_r$ and $N$ such that \begin{equation}\label{eq:beta_r}
    \min\{g(r-2),g'(r-2)+1\}\leq \beta_r\leq g(r-2),
    \end{equation}
    and for all $n\ge N$, we have 
    \[\sat(n,\mathcal{R}(K_r))= f(r-2)\left(n-f(r-2)\right)+\beta_r.\]
\end{lemma}

\begin{proof}
For $n\geq f(r-2)+2$, the construction $\Gamma_{r,n}$ demonstrates that \[\sat(n,\mathcal{R}(K_r))\leq f(r-2)(n-f(r-2))+g(r-2)\leq f(r-2)n.\]
Hence, by \cref{cor:structuralstability} with $\e=1$, there is some integer $N'$ such that for $n\geq N'$ and any $n$-vertex $\R(K_r)$-saturated edge-colored graph $\G$ with $\sat(n,\R(K_r))$-edges, there is a set $S\subseteq V(\G)$ of size $f(r-2)$ such that $S$ is complete to $V(\G)\setminus S$, and $|E(\G[S])|\geq g'(r-2)$.
Thus, we already have \[\sat(n,\R(K_r))\geq f(r-2)(n-f(r-2))+g'(r-2).\]
To show the existence of $\beta_r$, we now show that there exists $M$ such that 
the function ${\hat{\beta}_r:[M,\infty)\cap \mathbb{N} \rightarrow \mathbb{N}}$ given by $\hat{\beta}_r(n)=\sat(n,\R(K_r))-f(r-2)(n-f(r-2))$ is non-decreasing, and is therefore eventually equal to some constant $\beta_r$ between $g'(r-2)$ and $g(r-2)$.

As before, consider an $n$-vertex $\R(K_r)$-saturated edge-colored graph $\G$ with $\sat(n,\R(K_r))$ edges, where $n$ is large enough to deduce the existence of the set $S$ as above.
Note that there are at least $n-|S|-2g(r-2)$ vertices $v$ such that $N(v)=S$.
For each non-edge $e$ in $\binom{S}{2}$ and each color $c$, let $X_{c,e}$ be the vertex set of a copy of $\R(K_r)$ which is created when $e$ is added with color $c$. Note that for an arbitrary fixed color $c'$, we may choose $X_{c,e}=X_{{c'},e}$ for every color $c$ which does not appear in $\G[X_{{c'},e}]$.
Thus we may assume that the union $\mathcal{X}$ of all sets $X_{c,e}$ over all colors $c$ and all $e\in \binom{S}{2}\setminus E(G)$ has size at most $r\binom{r}{2}\binom{f(r-2)}{2}$.
In particular, if $n$ is large enough, there is some vertex $v\in V(\G)\setminus \mathcal{X}$ such that $N(v)=S$.
Now since $\G$ is $\R(K_r)$-saturated, so is $\G-v$, so $\sat(n-1,\R(K_r))\leq \sat(n,\R(K_r))-f(r-2)$, as promised.
This shows the existence of $\beta_r$ with $g'(r-2)\le \beta_r \le g(r-2)$. 

Finally, to show \eqref{eq:beta_r}, under the assumption that \[n>f(r-2)+2g(r-2)+2\left(f(r-2)\binom{f(r-2)}{2}\right)^{f(r-2)},\] suppose for contradiction that \[|E(\G)|= f(r-2)(n-f(r-2))+g'(r-2)<f(r-2)(n-f(r-2))+g(r-2).\]
Note that there is an independent set $X$ of size $2(f(r-2)\binom{f(r-2)}{2})^{f(r-2)}$ such that $N(v)=S$ for every $v\in X$.
Hence, by \cref{lem:rainbowlower}, $\G[S]\in \mathcal{F}_{r-2}$. Now since \[|E(\G[S])|\leq |E(\G)|-f(r-2)(n-f(r-2))=g'(r-2),\]
we have that $\G[S]\in \hat{\mathcal{F}}_{r-2}$ and so $|E(\G[S])|=g'(r-2)$ by the definition of $g'(r-2)$, and $V(\G)\setminus S$ is an independent set.
Since $g'(r-2)<g(r-2)$, it follows that $\G[S]$ is not $\R(K_{r-1})$-saturated.
Thus we can add some edge in some color to $S$ without creating $\R(K_{r-1})$ in $\G[S]$, and hence without creating $\R(K_r)$ in $\G$, a contradiction.
Therefore, this shows that the constant $\beta_r$ for which $\sat(n,\R(K_r))=f(r-2)(n-f(r-2))+\beta_r$ is between $\min\{g(r-2),g'(r-2)+1\}$ and $g(r-2)$, completing the proof.
\end{proof}

\subsection{Upper bound constructions for {\boldmath $f(k)$} and calculation of small values}
This subsection focuses on constructions that establish upper bounds on $f(k)$. However, we start by proving the following lower bounds on $f(k)$, which are equivalent to \cref{remark}.
\begin{prop}\label{prop:remark}
For every $k\ge 1$, we have
\begin{enumerate}
    \item $f(k)\ge k+1$ and
    \item $f(k)\ge k+2$ for every $k\ge 3$.
\end{enumerate}
\end{prop}

\begin{proof}
Property~\ref{item:prop2} immediately implies $f(k) \geq k+1$. We now show that for every $k\ge 3$, we have $f(k)\ge k+2$. Suppose, for a contradiction, that there is some $k\ge 3$ such that $f(k)=k+1$. Thus, $\Lambda_{k}$ has $k+1$ vertices. 
Since $k\geq 3$, Property~\ref{item:prop2} implies that $\Lambda_k$ has no missing edges.
Hence by Property~\ref{item:prop1}, there must be two edges $e,e'\in E(\Lambda_k)$ with the same color $c$. We split into two cases.
\begin{itemize}[leftmargin=*]
    \item If $e$ and $e'$ are not incident, then every copy of $\R(K_k)$ in $\Lambda_k$ uses either $e$ or $e'$, contradicting Property~\ref{item:prop3}.
    \item If $e$ and $e'$ are incident, then there must be a vertex $v$ that is adjacent to neither $e$ nor $e'$ (since $k+1\ge 4$). The graph $\Lambda_k-v$ contains both $e$ and $e'$ and thus is $\R(K_k)$-free, contradicting Property~\ref{item:prop2}.   
\end{itemize}
In both cases, we reached a contradiction, completing the proof of \cref{prop:remark}.
\end{proof}

The constructions for $f(1)$ and $f(2)$ are very simple.
\begin{lemma}\label{lem:F1F2}
    $f(1)=2$, $g'(1)=g(1)=0$,  $f(2)=3$, and $g'(2)=g(2)=3$.
\end{lemma}
\begin{proof}
Proposition~\ref{prop:remark} together with the construction $\Lambda_1=\overline{K_2}$ show that $f(1)=2$ and $g'(1)=g(1)=0$.
Proposition~\ref{prop:remark} together with the construction of a triangle with exactly two blue edges and one differently colored edge shows that $f(2)= 3$.
Furthermore, if $\Lambda$ is an edge-colored graph in $\hat{\mathcal{F}}_2$ with $V(\Lambda)=\{u,v,w\}$ and $uv\notin E(\Lambda)$, then $v$ violates Property~\ref{item:prop1}, a contradiction. 
Hence $g'(2)=g(2)=3$.
\end{proof}

We next give a recursive construction for general $k$.
\begin{lemma}\label{lem:Fk+1}
For each positive integer $k$, we have $\Gamma_{k+2,f(k)+2}\in \hat{\mathcal{F}}_{k+1}$.
In particular we have ${f(k+1)\leq f(k)+2}$, and if equality holds then $g(k+1)\leq \binom{f(k+1)}{2}-1$.
\end{lemma}
\begin{proof}
Recall that $\G_{k+2,f(k)+2}$ is $\R(K_{k+2})$-saturated, so in particular Property~\ref{item:prop1} holds. 
Let $u$ and $v$ be the two vertices in $V(\G_{k+2,f(k)+2})\setminus V(\Lambda_{k})$.
We can extend any $\R(K_{k})$ in $\Lambda_{k}$ to a copy of $\R(K_{k+1})$ in $\G_{k+2,f(k)+2}$ by adding either $u$ (in which case $v$ and all colors of edges incident to $v$ are avoided) or $v$ (in which case $v$ and all colors of edges incident to $v$ are avoided).
Thus Properties~\ref{item:prop2} and~\ref{item:prop3} are satisfied. Thus $|V(\G_{k+2,f(k)+2})|=f(k)+2$ is an upper bound for $f(k+1)$, and if this bound is tight then $|E(\G_{k+2,f(k)+2})|\leq \binom{f(k)+2}{2}-1$ is an upper bound for $g(k+1)$.
\end{proof}
We remark that based on the lower bound on $f(k)$ which we obtain in the next subsection, we have $f(k+1)=f(k)+2$ for infinitely many values of $k$.
This also allows us to push \cref{lem:F1F2} one step further, as follows.

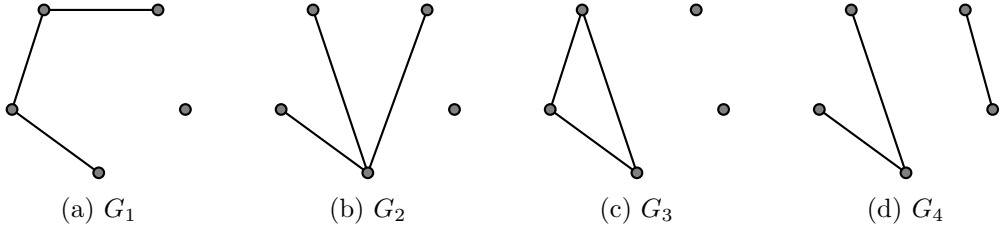
\begin{figure}[h]
\begin{subfigure}[b]{0.22\textwidth}
\centering
\begin{tikzpicture}[thick,scale=0.03]
\draw (50,10) node {} -- (12,38) node {} -- (26,82) node {} -- (76,82) node {} (88,38) node {} ;
\end{tikzpicture}
\caption{$G_1$}
\end{subfigure}
\begin{subfigure}[b]{0.22\textwidth}
\centering
\begin{tikzpicture}[thick,scale=0.03]
\draw (50,10) -- (12,38);
\draw (50,10) -- (26,82);
\draw (50,10) -- (76,82);
\draw (50,10) node {} (12,38) node {} (26,82) node {} (76,82) node {} (88,38) node {} ;
\end{tikzpicture}
\caption{$G_2$}
\end{subfigure}
\begin{subfigure}[b]{0.22\textwidth}
\centering
\begin{tikzpicture}[thick,scale=0.03]
\draw (50,10) -- (12,38);
\draw (50,10) -- (26,82);
\draw (12,38) -- (26,82);
\draw (50,10) node {} (12,38) node {} (26,82) node {} (76,82) node {} (88,38) node {} ;
\end{tikzpicture}
\caption{$G_3$}
\end{subfigure}
\begin{subfigure}[b]{0.22\textwidth}
\centering
\begin{tikzpicture}[thick,scale=0.03]
\draw (50,10) -- (12,38);
\draw (50,10) -- (26,82);
\draw (88,38) -- (76,82);
\draw (50,10) node {} (12,38) node {} (26,82) node {} (76,82) node {} (88,38) node {} ;
\end{tikzpicture}
\caption{$G_4$}
\end{subfigure}
\caption{Three edge graphs on five vertices}
\label{fig:3edgeGraphs}
\end{figure}

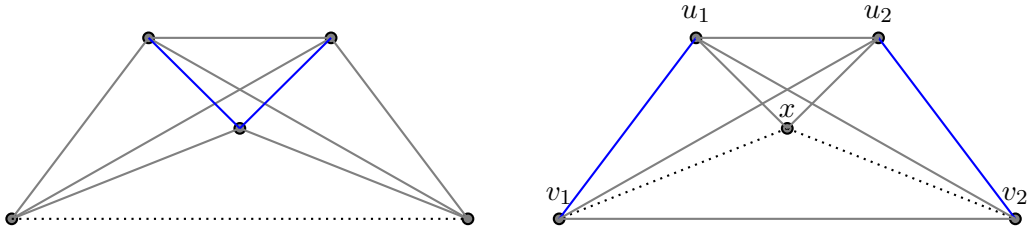
\begin{figure}[h]

\centering
\begin{tikzpicture}[thick,scale=0.12]
\draw (0,0) node {};
\draw (10,10) node {};
\draw (-10,10) node {};
\draw (-25,-10) node {};
\draw (25,-10) node {};
\draw [gray] (0,0) -- (-25,-10);
\draw [gray] (-25,-10) -- (-10,10);
\draw [gray] (0,0) -- (25,-10);
\draw [gray] (25,-10) -- (10,10);
\draw [gray] (10,10) -- (-10,10);
\draw [gray] (10,10) -- (-25,-10);
\draw [gray] (-10,10) -- (25,-10);
\draw [blue] (10,10) -- (0,0);
\draw [blue] (-10,10) -- (0,0);
\draw [dotted] (-25,-10) -- (25,-10);

\draw (60,0) node [label=$x$] {};
\draw (70,10) node [label=$u_2$] {};
\draw (50,10) node [label=$u_1$] {};
\draw (35,-10) node [label=$v_1$] {};
\draw (85,-10) node [label=$v_2$] {};
\draw [dotted] (60,0) -- (35,-10);
\draw [blue] (35,-10) -- (50,10);
\draw [dotted] (60,0) -- (85,-10);
\draw [blue] (85,-10) -- (70,10);
\draw [gray] (70,10) -- (50,10);
\draw [gray] (70,10) -- (35,-10);
\draw [gray] (50,10) -- (85,-10);
\draw [gray] (70,10) -- (60,0);
\draw [gray] (50,10) -- (60,0);
\draw [gray] (35,-10) -- (85,-10);
\end{tikzpicture}

\caption{The graphs $\Lambda_3$ (left) and $\Lambda'_3$ (right). The dotted lines represent non-edges, and each gray edge has a unique color.}
\label{fig:lambda_3}
\end{figure}

\begin{lemma}\label{lem:f3}
$f(3)=5$, $g(3)=9$, and $g'(3)=8$.
\end{lemma}
\begin{proof}
By \cref{prop:remark}, \cref{lem:F1F2}, and \cref{lem:Fk+1}, we have $f(3)=5$ and $g(3) \leq 9$.
Suppose for contradiction that there are two pairs $\{v,w\}$ and $\{v',w'\}$
of non-adjacent vertices in $\Lambda_3$, with $w'\notin \{v,w\}$.
Since $|V(\Lambda_3)|=f(3)=5$, there is some $u\in V(\Lambda_3)\setminus \{v,w,v',w'\}$.
Adding $vw$ to $\Lambda_3$ with the same color as $vu$ if $vu$ is an edge does not create a copy of $\R(K_4)$, contradicting the definition of $\Lambda_3$. If $vu$ is not an edge, adding $vw$ with any color does not create a copy of $\R(K_4)$, giving the same contradiction.
Thus, $g(3)=9$.

We now present a $5$-vertex graph in $\hat{\mathcal{F}}_3$ with exactly $8$ edges.
Let $V(\Lambda'_3) = \{x,v_1,v_2,u_1,u_2\}$.
The complement of $\Lambda'_3$ is exactly the edges $\{x,v_1\}$ and $\{x,v_2\}$.
The color of the edge $\{v_1,u_1\}$ is the same as the color $\{v_2,u_2\}$, and all other colors are distinct (see Figure~\ref{fig:lambda_3}).
It is straightforward to check that $\Lambda'_3 \in \hat{\mathcal{F}}_3$, and $|E(\Lambda'_3)|=8$.
Hence, $g'(3) \leq 8$.

To show that $g'(3) \geq 8$, let $\Lambda''$ be an edge-colored graph with $5$ vertices and $7$ edges; we will show that $\Lambda'' \notin \hat{\mathcal{F}}_3$.
There are $4$ distinct $3$ edge graphs on $5$ vertices; see Figure~\ref{fig:3edgeGraphs}.
Referring to the figure, if the complement of $\Lambda''$ is either $G_1$ or $G_3$, then every $K_3$ in $\Lambda''$ is incident to the right-most vertex; hence $\Lambda'' \notin \hat{\mathcal{F}}_3$.
If the complement of $\Lambda''$ is either $G_2$ or $G_4$, then the bottom vertex of $\Lambda''$ is not in any $K_3$. 
Hence, if $\Lambda''\in \hat{\mathcal{F}}_3$, then the graph obtained by deleting this vertex is also in $\hat{\mathcal{F}}_3$, contradicting the fact that $f(3)=5$.
Thus $\Lambda'' \notin \hat{\mathcal{F}}_3$, so $g'(3)=8$.
\end{proof}

\tikzstyle{every node}=[circle, draw, fill=black!50, inner sep=0pt, minimum width=4pt]

\tikzstyle{every node}=[circle, draw, fill=black!50, inner sep=0pt, minimum width=4pt]

To recap, the proof of \cref{lem:F1F2} gave an explicit construction for $\Lambda_2$: a triangle with two blue edges and one differently colored edge.
Thus we have an explicit construction of $\G_{2,5}$, which by \cref{lem:Fk+1} and \cref{lem:f3} is an explicit construction for $\Lambda_3$ (see Figure~\ref{fig:lambda_3}).
In turn, this gives us an explicit construction for $\G_{3,n}$.
In this case, we know that for sufficiently large $n$ this construction achieves the rainbow saturation number $\sat(n,K_5)$, based on \cref{thm:rainbow_construction} and \cref{lem:f3}.
It is tempting to conjecture that the graphs $\G_{r,n}$ are optimal for all $r\geq 3$ and $n$ sufficiently large (and thus that the constant $\beta_r$ in \cref{thm1} is always equal to $g(r-2)$).
To illustrate the difficulty in proving this, the graph $\Lambda'_3$ from \cref{lem:f3} (see Figure~\ref{fig:lambda_3}) can be used in an alternate construction for $\sat(n,\R(K_5))$ as follows. First, add vertices $z_1$ and $z_2$, complete to each other and $V(\Lambda'_3)$, such that the color of $z_1v_1$ is the same as that of $z_2v_2$,  the color of $z_1u_1$ is the same as that of $z_2u_2$ (but different to the color of $z_1v_1$), and all other new edges get their own unique color. Then add $n-7$ further vertices, each with neighborhood $V(\Lambda'_3)$, with all colors of edges incident to these vertices being unique in the entire construction.

We next give a construction that will be used to prove the upper bound in Theorem~\ref{thm1}. 
The construction is easiest to describe for $k=2 \binom{t}{2} + 1$ with $t \geq 3$.
Let $\Gamma'_k$ be a graph obtained by subdividing each edge of $\R(K_t)$ twice.
Note that $|V(\Gamma'_k)| = 2\binom{t}{2}+t$.
Each edge in $\Gamma'_k$ that is incident to a vertex of $R$ inherits the color of the subdivided edge of $R$, and each edge between subdivision vertices gets a unique color.
Add edges to $\Gamma'_k$ to make a complete graph $\Gamma_k$, giving each added edge a unique color.
It is straightforward to verify that $\Gamma_k \in \hat{\mathcal{F}}_k$, and a proof of this fact is given following the more general construction of \cref{lem:rainbowconstruction}.
Hence, $f(k) \leq 2\binom{t}{2} + t = k+(\sqrt{4k-3} - 1)/2$.

The following lemma generalizes the above construction for arbitrary values of $k \geq 4$.
\begin{lemma}\label{lem:rainbowconstruction}
    For any integer $k \geq 4$, we have
    \[f(k) \leq k+\left\lceil\frac{-1+\sqrt{4k-3}}{2} \right\rceil. \]
\end{lemma}

\begin{proof}
Let $n_k:=k+\left\lceil\frac{-1+\sqrt{4k-3}}{2} \right\rceil$.
Let $t:=\left\lceil \frac{1+\sqrt{4k-3}}{2}\right\rceil$. Note that  
\[{t \choose 2}\geq \frac{(\sqrt{4k-3}+1)(\sqrt{4k-3}-1)}{8}=\frac{k-1}{2}.\]
Since $k\geq 4$, we also have that $t\geq 3$ and \[{t \choose 2}\leq \left\lfloor \frac{(3+\sqrt{4k-3})(1+\sqrt{4k-3})}{8}\right\rfloor\leq \left\lfloor\frac{k+\sqrt{4k-3}}{2}\right\rfloor\leq \left\lfloor \frac{k+\sqrt{(k-0.1)^2}}{2}\right\rfloor=k-1.\]
Let $\ell:=2{t \choose 2}-(k-1)$ and let $m:=\frac{k-1-\ell}{2}$.
Let $\Gamma'_k$ be a subdivision of a copy of $K_t$, obtained by subdividing $m$ of its edges twice and subdividing the remaining $\ell$ edges once.
Let $X$ be the set of vertices in the original copy of $K_t$, and let $Y$ be the set of subdivision vertices.
Note that $|V(\Gamma'_k)|=t+(k-1)=n_k$. For each of the ${t \choose 2}$ paths $P$ of length $2$ or $3$ between vertices in $X$, color the first and last edges of $P$ with some color $c_P$ unique to $P$. Now add edges to make $\Gamma'_k$ complete and color all of the as yet uncolored edges with unique colors, so that the total number of colors used is ${{n_k}\choose {2}}-{t\choose 2}$.
Call this newly constructed edge-colored complete graph $\Gamma_k$. Now it is easy to verify that the minimum size of a hitting set for the paths of length two and three in $\Gamma'_k$ with repeated colors is $t-1$, and that the complement of any hitting set of this size induces a copy of $\R(K_k)$ in $\Gamma_k$. 
For any vertex $v\in V(\Gamma_k)$, let $x_1,x_2\in X$ be such that $v$ is in an $(x_1,x_2)$-path $P$ of length $2$ or $3$ in $\Gamma'_k$, and observe that $\Gamma_k-((X\setminus V(P))\cup \{v\})$ is a copy of $\R(K_k)$. This also shows that any color which appears on only one edge $e$ can be avoided, since we can take $v$ to be an endpoint of $e$. 
If $c$ is a repeated color, let $x$ be a vertex of $X$ which is not incident to any edge of color $c$, and observe that $\Gamma_k-(X\setminus \{x\})$ is a copy of $\R(K_k)$ avoiding $c$. This choice of $x$ is possible since $|X|=t\geq 3$.   
\end{proof}

As a simple corollary of \cref{prop:remark} and Lemmas~\ref{lem:sat_reduction}, \ref{lem:F1F2}, \ref{lem:Fk+1}, \ref{lem:f3}, and~\ref{lem:rainbowconstruction}, we can determine $\sat(n,\R(K_r))$ for small values of $r$ as follows.

\begin{theorem}\label{thm:smallvalues}
For sufficiently large $n$, the following holds.
\begin{equation*}
    \sat(n,\R(K_r))=\begin{cases}
    2n-4, & \text{if $r=3$},\\
    3n-6, & \text{if $r=4$},\\
    5n-16, & \text{if $r=5$, and}\\
    r(n-r)+\beta_r, & \text{if $r\in \{6,7,8,9\}$, with $\beta_r$ as in \cref{thm1}}.
    \end{cases}
\end{equation*}
\end{theorem}

\subsection{Lower bound of {\boldmath $f(k)$}}\label{sec:LowerBoundOffk}

In this subsection, we prove the following lemma which shows the lower bound in \cref{thm1}. 
\begin{lemma} \label{lem1}
For each positive integer $k$, we have
\[f(k) \geq k + (1/2)k^{1/3} - 1/2.\]
\end{lemma}

\begin{proof}
Define the \textit{rainbow-clique number} of an edge-colored graph $\G$ as the number of vertices in a largest rainbow-clique in $\G$.
Let $\G$ be an edge-colored graph of minimum order such that the rainbow-clique number of $\G$ is $k$ and, for every vertex $v \in V(\G)$, the rainbow-clique number of $\G-v$ is still $k$.
Note that $f(k) \geq |V(\G)| > k$.
Denote $|V(\G)| = n$; we show that $n \geq k + (1/2)k^{1/3} - 1/2$.

For any color $c$, let $E_c \subset E(\G)$ be the set of edges having color $c$.
Let $\mathcal{C}_2$ be the set of sets of edges with duplicate colors, together with pairs of vertices that are not edges in $\G$:
\[\mathcal{C}_2 = \{E_c: |E_c| \geq 2\} \cup \{\{\{u,v\}\}: \{u,v\} \notin E(\G) \}. \]
With abuse of terminology, we refer to any set $H$ such that the edge-colored graph induced on $V(\G) \setminus H$ is a $\R(K_k)$ as a {\em hitting set}.

For any set $S \subseteq V(\G)$ of vertices, let
\[I(S) = \sum_{E \in \mathcal{C}_2} \sum_{e \in E} \sum_{v \in S} \mathbbm{1}(v \in e), \]
where $\mathbbm{1}(v\in e)$ is the indicator function for the event that $v$ is incident to $e$, be the number of incidences between vertices of $S$ and the pairs counted in $\mathcal{C}_2$.
Note that $I(H)$ and $I(H')$ need not be equal for distinct hitting sets $H$ and $H'$.

The general plan is to give upper and lower bounds on $I(H)$, where $H$ is an arbitrary hitting set, as follows:
\begin{equation}\label{eq:hittingSetIncidenceBounds}
2(n-k)^3 + 3(n-k)^2 \geq I(H) \geq (1/4)n.
\end{equation}

It is straightforward to check that \cref{eq:hittingSetIncidenceBounds} is not satisfied for $k\leq n\leq k+(1/2)k^{1/3}-1/2$.
Thus, it remains to prove the upper and lower bounds of \cref{eq:hittingSetIncidenceBounds}.

Here comes the lower bound of \cref{eq:hittingSetIncidenceBounds}.
For each color $c$ that is used more than once, $H$ is incident to at least $|E_c|-1$ edges of $E_c$, since otherwise the graph induced on $V(\G) \setminus H$ will have two edges of color $c$.
Hence, if there are at least two edges of color $c$,
\[ \sum_{e \in E_c} \sum_{v\in H} \mathbbm{1}(v \in e) \geq |E_c|-1 \geq |E_c|/2.\]
Similarly, if $\{u,v\} \notin E(\G)$, then $H$ contains at least one of $u$ and $v$.

Since the rainbow-clique number of $\G$ is not changed when we remove a vertex, each vertex of $\G$ is incident to an edge with a duplicated color, or to a missing edge.
Hence,
\[|V| \leq \sum_{v \in V} \sum_{E \in \mathcal{C}_2} \sum_{e \in E} \mathbbm{1}(v \in e) = \sum_{E \in \mathcal{C}_2} 2|E|. \]
Combining these observations,
\[I(H) = \sum_{E \in \mathcal{C}_2} \sum_{e \in E} \sum_{v \in H} \mathbbm{1}(v \in e) \geq \sum_{E \in \mathcal{C}_2} |E|/2 \geq (1/4)|V|,  \]
which is the lower bound of \cref{eq:hittingSetIncidenceBounds}.

For the upper bound of \cref{eq:hittingSetIncidenceBounds}, we show that, for any fixed vertex $v \in V(\G)$,
\begin{equation}\label{eq:fixedVertex}
I(\{v\}) = \sum_{E \in \mathcal{C}_2} \sum_{e \in E} \mathbbm{1}(v \in e) \leq 3(n-k) + 2(n-k)^2.\end{equation}
Since $\G$ has minimum order, each vertex must be in some $\R(K_k)$. By definition, $|H| = n-k$ for each hitting set $H$.
Since \cref{eq:fixedVertex} holds for each vertex, taking the sum over the $n-k$ vertices in an arbitrary hitting set yields the upper bound of \cref{eq:hittingSetIncidenceBounds}.

Let 
\begin{align*}
\mathcal{C}_{v,1} &= \{E_c \in \mathcal{C}_2: \sum_{e \in E_c} \mathbbm{1}(v \in e) \geq 2\} \cup \{\{\{u,v\}\}: \{u,v\} \notin E(V)\}, \text{ and} \\
\mathcal{C}_{v,2} &= \{E_c \in \mathcal{C}_2:\abs{E_c}\ge 2, \sum_{e \in E_c} \mathbbm{1}(v \in e) = 1\}.\end{align*}
Note that 
\begin{equation}\label{eq:splitIv}
I(\{v\}) = \sum_{E \in \mathcal{C}_{v,1}} \sum_{e \in E} \mathbbm{1}(v \in e) + \sum_{E_c \in \mathcal{C}_{v,2}} \sum_{e \in E_c} \mathbbm{1}(v \in e).
\end{equation}
We separately bound each of the terms on the right side of \cref{eq:splitIv}.

If $\{u,v\} \notin E(\G)$ and $v \notin H$ for a hitting set $H$, then $u \in H$.
Let $E_{v,c} \subseteq E_c$ be those edges with color $c$ that are incident to $v$.
Since the graph induced on $V(\G) \setminus H$ can contain at most one edge of color $c$, all but one of the neighbors of $v$ via edges in $E_{v,c}$ must be in $H$.
Combining these observations, we have
\[\sum_{E \in \mathcal{C}_{v,1}} \sum_{e \in E} \mathbbm{1}(v \in e) \leq 2|H| = 2(n-k). \]

Assume that, for each color $c$ with $|E_c|\geq 2$, and each edge $e \in E_c$, it is not possible to give $e$ a new color (not present anywhere else in the graph) without increasing the rainbow-clique number of $\G$.
If this assumption does not hold for some edge $e$, then recolor $e$ to a new color.
This does not increase the order of $\G$, so it is sufficient to bound the order of the modified graph.
Since this recoloring strictly increases the number of colors and the number of colors cannot exceed the number of edges, the procedure will terminate.

It follows from this assumption that, for any $E_c \in \mathcal{C}_{v,2}$, there is an edge $e_c = \{x,y\} \in E_c$ with $x,y \neq v$ and $v$ is not adjacent to $y$ through an edge of color $c$, and a hitting set $H_c$ such that $v,y \notin H_c$.
Indeed, the assumption implies that there is a $K_{k+1}$ that includes $\{v,w\} \in E_c$ and is rainbow except for having two edges of color $c$.
We take $e_c$ to be the other edge of color $c$ in this $K_{k+1}$, and $H_c$ to be the complement of this $K_{k+1}$ together with $w$.

Fix a hitting set $H$ with $v \notin H$; this is possible since $v$ is contained in some $\R(K_k)$.
Partition $\mathcal{C}_{v,2}$ into $\mathcal{A}$ and $\mathcal{B}$, as follows.
If $E_c \in \mathcal{C}_{v,2}$ and there is a vertex of $H$ adjacent to $v$ through an edge of color $c$, then put $E_c \in \mathcal{A}$.
Otherwise, put $E_c \in \mathcal{B}$.
For $E_c \in \mathcal{B}$, $H$ does not contain a vertex incident to the edge of color $c$ that is incident to $v$.
Hence, $H$ must contain a vertex incident to $e_c$.

For each set $E_c \in \mathcal{A}$, there is a vertex $w_c \in H$ such that $\{v,w_c\} \in E_c$, and $w_c \neq w_{c'}$ for distinct colors $c \neq c'$.
Hence, $|\mathcal{A}| \leq |H| = n-k$.

For each set $E_c \in \mathcal{B}$, there is a vertex $y \in H$ such that $y \in e_c$.
For any vertex $y \neq v$, let $B_y = \{e_c: E_c \in \mathcal{B} \text{ and } y \in e_c\}$.
Since $\{e_c: E_c \in \mathcal{B}\} = \bigcup_{y \in H}B_y$ and $|\mathcal{B}|=|\{e_c: E_c \in \mathcal{B}\}|$, we have $|\mathcal{B}| \leq |H|\cdot \max_{y}|B_y| = (n-k)\cdot \max_y|B_y|$.

Fix a vertex $y \in H$ and a color $c$ such that $e_c \in B_y$.
By the definition of $B_y$, we have $y \in e_c$ and $v,y \notin H_c$.
Note that $H_c$ is incident to at least $|E_{c'}|-1$ edges of color $c'$ for each color $c'$.
Hence, for each edge $e_{c'} \in B_y$, it follows that $H_c$ contains at least one of the vertex $x_{c'}$ where $e_{c'} = \{y,x_{c'}\}$ and the vertex $w_{c'}$ where $\{v,w_{c'}\}$ is the edge of color $c'$ incident to $v$.
Furthermore, each vertex in $H_c$ is adjacent to $v$ by at most $1$ edge and adjacent to $y$ by at most one edge.
Combining these observations,
\[|B_y| \leq \sum_{E_{c'} \in B_y} \mathbbm{1}(x_{c'} \in H_{c}) + \sum_{E_{c'} \in B_y} \mathbbm{1}(w_{c'} \in H_{c}) \leq 2|H_c| = 2(n-k). \]
Hence, $|\mathcal{B}| \leq 2(n-k)^2$.

Taken together,
\[\sum_{E_c \in \mathcal{C}_{v,2}} \sum_{e \in E_c} \mathbbm{1}(v \in e) = |\mathcal{C}_{v,2}| = |\mathcal{A}| + |\mathcal{B}| \leq n-k+2(n-k)^2. \]
Combining this with \cref{eq:splitIv} and the previously obtained bound on $\sum_{E \in \mathcal{C}_{v,1}} \sum_{e \in E} \mathbbm{1}(v \in e)$, we have
\[I(\{v\}) \leq 3|H| + 2|H|^2 = 3(n-k) + 2(n-k)^2. \]
This finishes the proof of \cref{lem1}.
\end{proof}

Improving the preceding proof to give a bound on $\mathcal{B}$ of the form $|\mathcal{B}| \leq O(n)$ would lead to a bound of the form $f(k) \geq k + O(k^{1/2})$, which would determine the value of $f(k) - k$ up to a constant factor.

\subsection{Non-stability for rainbow saturation}
In this subsection, we aim to show that there is an abundance of constructions for $\R(K_r)$-saturated graphs. To motivate this, we start by mentioning the following result for ordinary saturation which shows that there is a linear gap between the smallest and the second smallest number of edges in an $n$-vertex $K_r$-saturated graph. 

\begin{theorem}[\cite{AFGS13}]
For every $r\ge 3$, if $G$ is an $n$-vertex $K_r$-saturated graph other than $K_{r-2}+\overline{K_{n-r+2}}$, then $|E(G)|\ge (r-1)n - \binom{r}{2} - 2$.
\end{theorem}

Contrary to the above, for sufficiently large $n$, there is an $n$-vertex $\R(K_r)$-saturated construction for any given number of edges between the saturation number and $\binom{n}{2}$.

\tikzstyle{every node}=[circle, draw, fill=black!50, inner sep=0pt, minimum width=4pt]

\begin{figure}[h]

\centering
\begin{tikzpicture}[thick,scale=0.15]
\draw [red] (0,0) -- (-25,-10);
\draw [red] (-25,-10) -- (-10,10);
\draw [blue] (0,0) -- (25,-10);
\draw [blue] (25,-10) -- (10,10);
\draw [green] (10,10) -- (-10,10);
\draw [green] (10,10) -- (-25,-10);
\draw [green] (-10,10) -- (25,-10);
\draw [green] (10,10) -- (0,0);
\draw [green] (-10,10) -- (0,0);
\draw [dotted] (-25,-10) -- (25,-10);
\draw (0,0) node {};
\draw (10,10) node {};
\draw (-10,10) node {};
\draw (-25,-10) node {};
\draw (25,-10) node {};
\end{tikzpicture}

\caption{A $\R(K_3)$-saturated graph with exactly one missing (dotted) edge, denoted $\Lambda'_3$}
\label{fig:nice_small}
\end{figure}
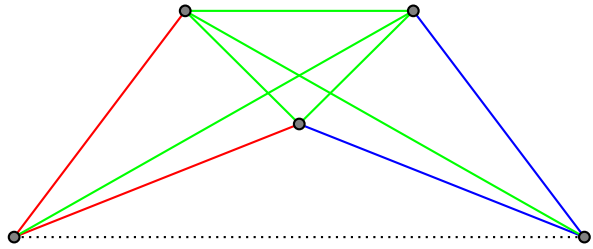

\begin{theorem}\label{thm:non_stability}
    For every $r\geq 3$, there exists $N'$ such that for all $n\geq N'$ and every integer $m$ satisfying $\sat(n,\R(K_r))\leq m\leq \binom{n}{2}$, there is a $\R(K_r)$-saturated graph with exactly $n$ vertices and $m$ edges.
\end{theorem}
\begin{proof}
First choose $N\geq f(r-2)+2\beta_r+5$ to be sufficiently large in the sense of both \cref{cor:structuralstability} (with $\e=1$) and \cref{lem:sat_reduction}, and let $E_r:=\binom{N}{2}-f(r-2)(N-f(r-2))-\beta_r$.
Let $C^*:=4(N-f(r-2))+10$.
We will take $N':= \max\{N,(2r-1)(E_r+C^*(N-f(r-2))+\binom{C^*}{2})\}$.
Let $M=\binom{n}{2}-m$ (the target number of missing edges).
Note that if $n=N+C$ for $C\geq 0$, then 
\[M\leq \binom{n}{2}-f(r-2)(N+C-f(r-2))-\beta_r=E_r+C(N-f(r-2))+\binom{C}{2}.\]
We split into two cases depending on $M$. In each case, we leave it to the readers to verify that the constructions are, in fact, $\R(K_r)$-saturated with the desired number of edges.

{\bf Case 1:} $M\geq E_r+C^*(N-f(r-2))+\binom{C^*}{2}$. Let $a$ be the minimum positive integer such that $M\leq E_r+a(N-f(r-2))+\binom{a}{2}$, and let $b$ be such that $M=E_r+a(N-f(r-2)) + \binom{a}{2}-b$.
Note that $n\geq N+a$ by our previous observation.
Observe that $0\leq b\leq (N-f(r-2))+a-1$ by our choice of $a$, and $a\geq C^*=4(N-f(r-2))+10$ by assumption.
Now we can write $b$ as $6x+y$ for some non-negative integers $x$ and $y$ with $y\leq 5$.
Now if $x\geq N-f(r-2)$, then $x\geq 5\geq y$ and $4x+2y\leq b-x\leq a$. Otherwise, $4x+2y\leq C^*\leq a$.

We now describe the construction for this case.
Let $n':=N+a-3x-y$ and let $G'$ be an edge-colored graph with $n'$ vertices and $\sat(n',\R(K_r))$ edges which is $\R(K_r)$ saturated.
There is a set $S\subset V(G')$ of size $f(r-2)$ which is complete to $V(G')\setminus S$, and at most $\beta_r$ edges of $G'$ have both endpoints in $V(G')\setminus S$. 
Since $N> f(r-2) + 2\beta_r$, there are disjoint subsets $A$, $X$, and $Y$ of $V(G')$ of sizes $a-4x-2y$, $x$, and $y$ respectively such that $A\cup X\cup Y$ is independent in $G'$ and $N(v)=S$ for all $v\in A\cup X\cup Y$.
For each vertex $v\in X$, we add three new vertices $u_{v,1},u_{v,2}$, and $u_{v,3}$, such that for each $i\in [3]$, we have $N[u_{v,i}]=N[v]\cup \{u_{v,1},u_{v,2},u_{v,3}\}$,
and the color of $u_{v,i}w$ is the color of $vw$ for all $w\in N_{G'}(v)$, and all edges in the graph induced on $\{v,u_{v,1},u_{v,2},u_{v,3}\}$ are red (an arbitrarily chosen color which may appear elsewhere in the construction).
For each vertex $v\in Y$, we add one new vertex $u'_{v}$, such that for each $N(u'_v)=N[v]$,
and the color of $u'_{v}w$ is the color of $vw$ for all $w\in N_{G'}(v)$, and the edge $vu'_v$ is red.
Since the graph $G'$ is $\mathcal{R}(K_r)$-saturated, so is the newly constructed graph. 
The number of missing edges is $E_r+a(N-f(r-2))+\binom{a}{2}-6x-y=M$.
We complete the construction by adding $n-N-a$ dominant vertices, so that all edges incident to them are red.

{\bf Case 2:} $M\leq E_r+C^*(N-f(r-2))+\binom{C^*}{2}$.
We first construct a $\R(K_r)$-saturated graph with at most $(2r-1)M$ vertices and $M$ missing edges as follows. 
If $r\geq 4$, let $H_r$ be the $2r-2$ vertex graph $\overline{K_2}+2K_{r-2}$. Observe that $\R(H_r)$ has a unique non-edge which cannot be added without creating $\R(K_r)$. Add all other non-edges in color red, and call this edge-colored graph $\Lambda'_r$. If $r=3$, instead let $\Lambda'_r$ be the edge-colored graph depicted in Figure~\ref{fig:nice_small}.
We take $M$ disjoint copies of $\Lambda'_r$, and make them complete to each other with red edges to obtain a $\R(K_r)$-saturated graph with exactly $M$ missing edges. We then add $n-|V(\Lambda'_r)|M$ dominant vertices so that all edges incident to them are red.
\end{proof}

\section{{\boldmath $(H,1)$-saturation and -semisaturation, with applications to rainbow saturation problems}}\label{sec:MoreProofs}

In this section, for each $r\geq 3$ and all sufficiently large $n$, we determine $\sat_1(n,K_r)$, $\ssat_1(n,K_r)$, $\sat(\R(K_n),\R(K_r))$, and $\ssat(n,\R(K_r))$.
The first ingredient we will need is \cref{prop:comparison_semi}, which we now prove.

\begin{proof}[Proof of \cref{prop:comparison_semi}]
Consider removing an edge $e$ from $G$ and then adding edges $e_1$ and $e_2$ to $G$. We aim to show that this creates a new copy of $H$ in $G$. Without loss of generality, we assume that $e_1\neq e$. Now, in $\R(G)$, we add the edge $e_1$ with the same color as $e$. This must create a new copy of $\R(H)$ and moreover this copy does not use the edge $e$. Thus, removing $e$ and adding $e_1$ in $G$ creates a new copy of $H$, which completes the proof.
If we additionally assume that $\R(G)$ is $\R(H)$-saturated, then $G$ is also $H$-free, and thus $(H,1)$-saturated.
\end{proof}

Our strategy to prove the results of this section is to apply \cref{lem:general} with an appropriately chosen $k$, $\ell$, and $\mathcal{F}$ to obtain a lower bound which matches the upper bounds given by our constructions. 
In order to do this, the following lemma is essential for determining the constant $k$ and the family $\mathcal{F}$.
In this way, it is analogous to Lemma~\ref{lem1} in the previous section.

\begin{lemma} \label{lem2}
    If $G$ is a graph such that $\omega(G-v)=\omega(G)$ for each for each vertex $v$ of $G$, then the complement $\overline{G}$ of $G$ contains a matching of size $\omega(G)$. In particular, $|V(G)|\geq 2\omega(G)$. 
\end{lemma}

\begin{proof}
Let $H$ be a clique of order $\omega(G)$ in $G$ and let $M=\{x_iy_i:i\in [t]\}$ be a matching in $\overline{G}$ of maximum possible size subject to the condition that $X:=\{x_i:i\in [t]\}\subseteq V(H)$ and $V(H)\setminus X$ is complete to $\{y_i:i\in [t]\}$ in $G$. Suppose for contradiction that $V(H)\setminus X$ contains some vertex $v$, and let $H'$ be clique of order $\omega(G)$ in $G-v$. 
Let $A:=V(H)\setminus V(H')$, let $B:=V(H')\setminus V(H)$, and note that $G':=\overline{G}[A\cup B]$ is bipartite with bipartition $(A,B)$.
If for some $S\subseteq A$ we have $|N_{G'}(S)|<|S|$, then $(V(H')\setminus N_{G'}(S))\cup S$ induces a clique in $G$ of order greater than $n$, a contradiction. Hence, by Hall's marriage theorem, there is a perfect matching $M'$ from $A$ to $B$ in $G'$.
Let $\{v_i:i\in [t']\}:=A\setminus X$, and for each $i\in [t']$ let $w_i$ be the vertex in $B$ which is matched to $v_i$ in $M'$.
Note that $t'\ge 1$, since $v\in A\setminus X$. Also,
since $V(H)\setminus X$ is complete to $\{y_i:i\in [t]\}$ in $G$, the set $\{w_i:i\in [t']\}$ is disjoint from $\{y_i:i\in [t]\}$, and so $M'\cup M$ is a matching in $\overline{G}$ which is strictly larger than $M$. 
Finally, observe that $V(H)\setminus (X\cup \{v_i:i\in [t']\})$ is complete to $\{y_i:i\in [t]\}$ in $G$ by the definition of $M$, and that $V(H)\setminus (X\cup \{v_i:i\in [t']\})$ is complete to $\{w_i:i\in [t']\}$ in $G$ since $H'$ is a clique containing both sets.
Hence $M\cup M'$ contradicts our choice of $M$.
Therefore $V(H)\setminus X$ is empty, and so $M$ is the desired matching of size $\omega(G)$ in $\overline{G}$.
\end{proof}

\subsection{{\boldmath $(K_r,1)$}-semisaturation and {\boldmath $\mathcal{R}(K_r)$}-semisaturation}

In this subsection, we prove \cref{thm:semi-saturation,thm2'} simultaneously.
For an integer $n\geq 3$, we define $G_{3,n}$ to be the complete bipartite graph $K_{2,n-2}$. For integers $r\geq 4$ and $n\geq r$, we define $G_{r,n}$ to be the complete $r$-partite graph $K_{r-1} +\overline{K_{n+1-r}}$.

\begin{theorem}\label{thm:stabsemisat}
For every integer $r\ge 3$ and real number $\epsilon\in (0,1]$, there exists $N\ge r+1$ such that for all $n\ge N$ and every $n$-vertex graph $G$ with at most $(r-\epsilon)n$ edges,
the following are equivalent:
\begin{enumerate}
    \item \label{item:1} $G$ is $(K_r,1)$-semisaturated;
    \item \label{item:2} $G_{r,n}\subseteq G$;
    \item \label{item:3} $\mathcal{R}(G)$ is $\mathcal{R}(K_r)$-semisaturated.
\end{enumerate}
In particular,
\begin{equation*}
    \ssat(n,\mathcal{R}(K_r))=\ssat_1(n,K_r)=\begin{cases}
    2(n-2), & \text{if $r=3$}.\\
    (r-1)(n-r+1)+\binom{r-1}{2}, & \text{otherwise}.
    \end{cases}
\end{equation*}

\end{theorem}

In view of the graph $\overline{K_{n-r}}+\overline{K_2}+K_{r-2}$, the above theorem does not hold with $\e=0$.
We define $\mathcal{F}'_3$ to be the class of all $2$-vertex graphs, and for $r\geq 4$, we define $\mathcal{F}'_r:=\{K_{r-1}\}$.
\begin{lemma}\label{lem:semi_lower}
    For $r\ge 3$, given an $n$-vertex $(K_r,1)$-semisaturated graph $G$,
    there is no set $S\subseteq V(G)$ of size $r-1$ for which there is an independent set $X$ in $G-S$ of size $r+1$ such that for some $v\in X$ and every $w\in X\setminus \{v\}$, we have that $N(v)\cap N(w)\subseteq S$ and $G[N(v)\cap N(w)]$ is not isomorphic to a graph in $\mathcal{F}'_r$.
\end{lemma}

\begin{proof} 
Suppose for contradiction that such an $S$, $X$, and $v$ exist.
Now for every $w\in X\setminus \{v\}$, we have that $\omega(G[N(v)\cap N(w)])\leq r-2$ and $|N(v)\cap N(w)|<2(r-2)$. Hence by \cref{lem2}, there is a vertex $s_w\in S$ such that there is no copy of $K_{r-2}$ in $G[(N(v)\cap N(w))\setminus \{s_w\}]$ (where $s_w$ can be chosen arbitrarily if $N(v)\cap N(w)=\emptyset$). 
Now by the pigeonhole principle, there are distinct vertices $w$ and $w'$ in $X$ such that $s_{w}=s_{w'}$.
By our choice of $s_w$ and the fact that $w$ and $w'$ are not adjacent, we can subtract the edge $vs_w$ from $G$ and then add the edges $vw$ and $vw'$ without creating a new copy of $K_r$.
This contradicts that $G$ is $(K_r,1)$-semisaturated, completing the proof.
\end{proof}

\begin{proof}[Proof of \cref{thm:stabsemisat}]
Pick $N$ according to \cref{lem:general} with $k:=r-1$ and $\ell:=r+1$.
Let $\mathcal{F}$ be the class of all $2$-vertex graphs if $r=3$, and otherwise $\mathcal{F}:=\{K_{r-1}\}$.
Now \cref{lem:general} and \cref{lem:semi_lower} together show that (\ref{item:1}) implies (\ref{item:2}).
Assuming (\ref{item:2}), adding any edge to $G$ creates $r-1$ new copies of $K_r$, such that no edge of $G$ is in their common intersection other than the newly added edge. Since each color appears at most once in $\mathcal{R}(G)$, it follows that one of these copies avoids the color of the newly added edge. Thus (\ref{item:2}) implies (\ref{item:3}).
Finally, it follows from \cref{prop:comparison_semi} that (\ref{item:3}) implies (\ref{item:1}).
\end{proof}

\subsection{{\boldmath $(K_r,1)$}-saturation and {\boldmath $(\R(K_n),\mathcal{R}(K_r))$}-saturation}

In this subsection, we prove \cref{thm:rainbow_construction,thm2} simultaneously.
Given positive integers $r\ge 3$ and $n\geq 2r-3$, we define $G'_{r,n}$ to be the complete join of $\overline{(r-2)K_2}$ with an independent set of size $n-2(r-2)$.
\begin{theorem}\label{thm:stabsat}
For every integer $r\ge 3$ and real number $\epsilon\in (0,1]$, there exists $N\ge 2(r-2)+2$ such that for all $n\ge N$ and every $n$-vertex graph $G$ with at most $(2(r-2)+1-\epsilon)n$ edges,
the following are equivalent:
\begin{enumerate}
    \item \label{item:1a} $G$ is $(K_r,1)$-saturated;
    \item \label{item:2a} $G\cong G'_{r,n}$;
    \item \label{item:3a} $\mathcal{R}(G)$ is $\mathcal{R}(K_r)$-saturated.
\end{enumerate}
In particular,
\begin{equation*}
    \sat(\R(K_n),\mathcal{R}(K_r))=\sat_1(n,K_r)=2(r-2)(n-r+1).
\end{equation*}

\end{theorem}

In view of the graph $\overline{K_{n+3-2r}}+\overline{K_3}+\overline{(r-3)K_2}$, the above theorem does not hold with $\e=0$. Given an integer $r\geq 3$, we define $\mathcal{F}''_r$ to be the family of $2(r-2)$-vertex graphs $G$ such that $\omega(G)=r-2$ and $\overline{G}$ has a perfect matching.

\begin{lemma}\label{lem:1satX}
Given positive integers $n\geq r\geq 3$ and an $n$-vertex $(K_r,1)$-saturated graph $G$, there is no set $S\subseteq V(G)$ of size $2(r-2)$ for which there is an independent set $X$ in $G-S$ of size $2r-2$ and a vertex $v\in X$, such that for every $w\in X\setminus \{v\}$ we have that $N(v)\cap N(w)\subseteq S$ and $G[N(v)\cap N(w)]$
is not isomorphic to a graph in $\mathcal{F}''_{r}$.
\end{lemma}

\begin{proof}
Suppose for contradiction that such an $S$, $X$, and $v$ exist. 
Consider an arbitrary vertex $w\in X\setminus \{v\}$.
Since $(K_r,1)$-saturation implies $K_r$-saturation, we have $\omega(G[N(v)\cap N(w)])= r-2$.
Hence, by Lemma~\ref{lem2} there is some $s_w\in S$ such that $\omega(G[N(v)\cap N(w)]-s_w)<r-2$ (where $s_w$ can be chosen arbitrarily if $N(v)\cap N(w)=\emptyset$). 
By the pigeonhole principle, there are distinct vertices $w$ and $w'$ in $X\setminus\{v\}$ such that $s_w=s_{w'}$.
Now, adding the edges $vw$ and $vw'$ to $G$ and subtracting the edge $vs_w$ does not create a copy of $K_r$, a contradiction.
\end{proof}

\begin{proof}[Proof of \cref{thm:stabsat}]
We pick $N$ according to \cref{lem:general} with $k:=2(r-2)$, $\ell:=2r-2$.

First suppose (\ref{item:1a}) holds. By \cref{lem:general} together with \cref{lem:1satX}, there is a subset $S$ of $G$ of size $2(r-2)$ such that $S$ is complete to $V(G)\setminus S$, and $G[S]\in \mathcal{F}''_{r}$. Since $G$ is $K_r$-free, it follows that there is no edge in $G-S$. Thus $G$ is a subgraph of $G'_{r,n}$, and since $G'_{r,n}$ is $K_r$-free and $G$ is $K_r$-saturated, we have $G\cong G'_{r,n}$. Thus (\ref{item:1a}) implies (\ref{item:2a}).

Since $N\geq 2(r-2)+2$, it is easy to check that adding any edge to $G'_{r,n}$ creates two copies of $K_r$ which share no edge other than the newly added edge. 
Thus adding any edge with any color to $\R(G'_{r,n})$ creates a copy of $\R(K_r)$, so (\ref{item:2a}) implies (\ref{item:3a}).

Finally, by \cref{prop:comparison_semi} we have that (\ref{item:3a}) implies (\ref{item:1a}).
\end{proof}

\section{Concluding remarks}\label{sec:Conclusion}
It would be interesting to close the gaps between the lower and upper bounds of $\sat(n,\R(K_r))$ obtained in Theorem~\ref{thm1}. By Lemma~\ref{lem:sat_reduction}, this is equivalent to finding better bounds on $f(k)$.

\begin{conjecture}
For every $k\geq 4$, $\alpha_{k+2}=f(k)=k+\left\lceil\frac{-1+\sqrt{4k-3}}{2}\right\rceil$.
\end{conjecture}
More strongly, we conjecture that, for $k$ of the form $2\binom{t}{2}+1$ where $t \geq 3$, the construction given in \cref{lem:rainbowconstruction} is the unique graph in $\hat{\mathcal{F}}_k$ with $f(k)$ vertices.

\cref{cor:comparison} shows a linear gap between $\sat(n,\R(K_r))$ and $\ssat(n,\R(K_r))$ for $r\ge 5$. By \cref{cor:comparison,thm2',ehm}, for every $r\ge 3$, we have \[(r-1)n+O(1)=\ssat(n,\R(K_r))\ge \wsat(n,\R(K_r))\ge \wsat(n,K_r)=(r-2)n+O(1).\] 
Thus it would be interesting to determine the gap between $\ssat(n,\R(K_r))$ and $\wsat(n,\R(K_r))$. For ordinary saturation, the known proofs to determine $\wsat(n,K_r)$ use algebraic techniques, see, e.g. \cite{A85,FFS11}, but it is unclear whether such methods can be generalized to the rainbow setting. 
In a similar vein to Question~6.2 of \cite{BJLMO22}, we ask the following.

\begin{question}
Given $r\geq 3$, is the following true?
\[\ssat(n,\R(K_r))= \wsat(n,\R(K_r))+o(n).\]
\end{question}

It is possible to generalize the notion of $(H,1)$-saturation in the following manner. For any $k\ge 0$, a graph $G$ is called $(H,k)$-saturated if $G$ is $H$-free, but removing any $k$ edges and then adding any $k+1$ edges (possibly including removed edges) creates a copy of $H$. Thus, ordinary $H$-saturation is the same as $(H,0)$-saturation. Let $\sat_k(n,H)$ denote the minimum number of edges in an $n$-vertex $(H,k)$-saturated graph. Note that for every $0\le l\le k$, if $G$ is $(H,k)$-saturated, then it is $(H,l)$-saturated; thus, $\sat_l(n,H)\le \sat_k(n,H)$. The upper bound construction for $(K_r,1)$-saturation can be extended to obtain an upper bound for $\sat_k(n,K_r)$.

\begin{prop} 
For every $r\ge 3$, $k\ge 0$, and $n\ge (k+1)(r-2)$, the saturation number $\sat_k(n,K_r)\le (k+1)(r-2)\left(n - (k+1)(r-2)\right) + (k+1)^2\binom{r-2}{2}$. 
\end{prop}

The above upper bound is achieved by the construction $\overline{K_{n-(k+1)(r-2)}}+ \overline{(r-2)K_{k+1}}$. 
However, the lower bound argument for Theorem~\ref{thm2} does not seem to generalize to arbitrary~$k$.
A natural generalization of Lemma~\ref{lem2} for this purpose would be that, if $G$ is a graph such that $\omega(G) = \omega(G-S)$ for every $S \subseteq V(G)$ of size at most $t$, then $|V(G)| \geq (t+1)\omega(G)$.
However, this is not true. 
For example, if $G$ is the complement of the Petersen graph, then $\omega(G-S)=4$ for every $S \subset V(G)$ with $|S|\leq 2$, and $|V(G)|=10<12$.

\begin{question}
Given integers $r\ge 2$ and $t\ge 2$, what is the minimum integer $n$ such that some $n$-vertex graph $G$ satisfies $\omega(G)=\omega(G-S)=r$ for every $S\subseteq V(G)$ of size at most $t$?
\end{question}

It is known that the ordinary saturation numbers $\sat(n,H)$ are always at most linear in~$n$, see, e.g., \cite{FFS11,KT86}. It would be thus interesting to investigate whether the other variants of saturation numbers also behave in the same way. Indeed, Behague, Johnston, Letzter, Morrison, and Ogden proved the following (which was originally conjectured by Gir\~{a}o, Lewis, and Popielarz).

\begin{theorem}[\cite{BJLMO22,GLP20}]
For every graph $H$, we have $\sat(n,\R(H))=O(n)$.
\end{theorem} 

Remember that $\sat(\R(K_n),\R(H))$ can be infinite (for example, when $H=K_{1,4}$ and odd $n\geq 5$). 
However we suspect that $\sat(\R(K_n),\R(H))$ is finite for infinitely many $H$.
Indeed, it may be true that for every $H$ there exists a constant $C(H)$ such that ${\sat(\R(K_n),\R(H))\leq C(H) n}$ for infinitely many $n$ (note that $\sat(\R(K_n),\R(K_{1,4}))= 3n/2$ for all even $n\geq 4$). 
Instead of considering this problem directly, we may instead consider the closely related $1$-saturation number $\sat_1(n,H)$, which has the advantage of always being finite (indeed it is bounded above by the extremal number $\ex(n,H)$).
Here, we ask the following.
\begin{question}\label{qst:sat1}
Given a graph $H$, does there exist a constant $c(H)$ such that $\sat_1(n,H)\leq c(H)n$ for infinitely many $n$ (or indeed for all positive integers $n$)?
\end{question}

It is natural to consider the analogous question for saturation numbers $\sat_k(n,H)$ with $k\geq 2$. However we make the following preliminary observation, answering this question in the negative for $k\geq 3$. 

\begin{theorem} \label{sat_c4}
The saturation number $\sat_3(n,C_4)$ is superlinear. In particular, we have \[\sat_3(n,C_4)= \Omega(n^{20/19}).\]
\end{theorem}
\begin{proof}
Let $G$ be an $n$-vertex graph which is $(C_4,3)$-saturated. 
Given a pair of non-adjacent vertices $v$ and $w$ in $G$, we say that $\{v,w\}$ is guarded by $xy\in E(G)$ if $xy$ is the central edge of a path of length $3$ from $v$ to $w$ in $G$ (note that such a path must exist in a $C_4$-saturated graph).
Let $v$ be a vertex of maximum degree in $G$, and consider $G_v:=G[N(v)]$. Since $G$ has no $C_4$, every component of $G_v$ has size at most $2$. 

\begin{claim}\label{clm:externalgaurds}
For every four distinct components $X_1,X_2,X_3$, and $X_4$ of $G_v$ there exist distinct $i,j\in [4]$ and $x\in X_i, y\in X_j$ such that the pair $\{x,y\}$ is guarded by an edge in $E(G-N[v])$.
\end{claim}
\begin{proof}
Since $G$ is $C_4$-saturated, for every distinct $i,j\in [4]$ and $x\in X_i, y\in X_j$, we have that $\{x,y\}$ is guarded by an edge $e$ in $E(G)$. 
If $|X_i|=|X_j|=1$, then there is no such path of length $3$ in $G$ containing $x,y$, and $v$. Let $P$ be an $(x,y)$-path of length~$3$ in~$G$. Since every neighbor of $x$ or $y$ in $G-v$ is in $V(G)\setminus N[v]$, the edge of $P$ which guards $\{x,y\}$ is in $E(G-N[v])$ as required. 
Hence, we may assume $|X_1|=2$. 
Let $X_1=\{a_1,b_1\}$ and for each $i\in [4]$, let $X_i$ contains a vertex $a_i$. Consider the graph $G'$ with $V(G'):=V(G)$ and $E(G'):=(E(G)\setminus\{va_1,vb_1,va_2\})\cup \{a_1a_2,b_1a_2,a_1a_3,b_1a_4\}$.
Now $G'$ contains a subgraph $H\cong C_4$, since $G$ is $(C_4,3)$-saturated. However, it is quick to check that $C_4$ is not a subgraph of $G'[N_G[v]]$. 
Hence, $H$ contains some vertex $w\in V(G)\setminus N_G[v]$ and some edge $xy\in \{a_1a_2,b_1a_2,a_1a_3,b_1a_4\}$. Since $G$ has no $C_4$ subgraph, $w$ has at most one neighbor in $N_G(v)$, so some neighbor $w'$ of $w$ in $H$ is in $V(G)\setminus N_G[v]$. 
Now $ww'$ guards $\{x,y\}$, as required.
\end{proof}

\begin{claim}\label{maxdegreebound}
The maximum degree $\Delta$ of $G$ satisfies $\Delta<\sqrt{24|E(G)|}$
\end{claim}

\begin{proof}
Let $S$ be the set of pairs of vertices $\{x,y\}\subseteq N(v)$ such that $xy\notin E(G)$ and $\{x,y\}$ is guarded by an edge in $E(G-N[v])$. 
Note that every edge $uw \in E(G-N[v])$ guards at most one pair in $S$, since otherwise either $u$ or $w$ is adjacent to two neighbors of $v$, contradicting the fact the $G$ does not contain $C_4$.

Hence $|S|\leq |E(G-N[v])|$. 
Let $G'$ be the graph whose vertices are the components of $G_v$ where two vertices $X_i$ and $X_j$ are adjacent whenever there exist $x\in X_i$ and $y\in X_j$ with $\{x,y\}\in S$. 
Now $|V(G')|\geq \Delta/2$ and $|S|\geq |E(G')|$. By Claim~\ref{clm:externalgaurds} and Turan's Theorem~\cite{T41} for $K_4$, we have
\[|E(G)|\geq \Delta+ |E(G-N[v])| \geq \Delta +|E(G')|\geq \Delta+ 3{{\lfloor |V(G')|/3\rfloor}\choose {2}}>\frac{\Delta^2}{24}.\qedhere\]
\end{proof}

First, $\Delta \le \frac{n^{10/19}}{2}$ because otherwise \cref{sat_c4} follows from Claim~\ref{maxdegreebound}. 
Now, let $C$ be the set of vertices of degree at least $n^{8/19}$ in $G$. Clearly, $|C|\le n^{12/19}$ because otherwise the number of edges in $G$ is at least $\frac{n^{20/19}}{2}$. Let $A=E(G[C])$. Then, we have $|A|\le n^{18/19}$ by the known result on the extremal number of $C_4$ (i.e., $\ex(n,C_4)\leq \frac{1}{2}n^{3/2} + o(n^{3/2})$, see \cite{FS13,KST54} for references).

Let $B=E(G)\setminus A$.
Now every pair $\{x,y\}$ in ${{V(G)}\choose {2}}\setminus E(G)$ is guarded by some edge in $E(G)$, and every edge $uw\in E(G)$ guards at most $\deg(u)\deg(w)$ edges in ${{V(G)}\choose {2}}\setminus E(G)$.
Hence we obtain by the bounds on $\Delta$ and $|A|$
\[{{n}\choose {2}}-|A|-|B|\leq |B|\Delta n^{8/19}+|A|\Delta^2\le |B|\frac{n^{18/19}}{2} + \frac{n^2}{4}.\]
Hence, we obtain $|B| \ge \frac{n^{20/19}}{2}$, thus the number of edges in $G$ is $\Omega(n^{20/19})$.
\end{proof}

We finish by suggesting a couple of other directions for future research. There are studies of generalized rainbow Tur\'an problems, see, e.g., \cite{GMMP22,J22}. There are also studies on generalized saturation problems (where the objective is to minimize the number of copies of a fixed graph instead of edges) see, e.g., \cite{CL20,KMTT20}. It might be worth studying the generalized saturation problems in rainbow settings.

In another direction, graph saturation has been extended to the setting where the addition of the edges is restricted to some host graph other than the complete graph $K_n$. There are studies with several several different host graphs such as complete bipartite graphs \cite{CCH21,GKS15}, complete multipartite graphs \cite{FJPW16,GKP19,R17}, hypercubes \cite{CG08,JP17,MNS17} and random graphs \cite{DSZ22,KS17}. Rainbow saturation problems can be studied on such host graphs.

\section*{Acknowledgment}
We are grateful to the Discrete Mathematics Group of the Institute for Basic Science for organizing the 2021 DIMAG Internal Workshop at Gangneung, where the authors initiated this project, and to the participants of this workshop for stimulating discussions. We are also thankful to the anonymous referee for their careful reading, and in particular for pointing out a small error in the proof of \cref{thm:non_stability} in a previous version of this paper and observing that our argument for \cref{lem1} yields a better second-order term than we had previously claimed.

\bibliographystyle{abbrv}
\bibliography{references.bib}

\end{document}